\documentclass[10pt]{amsart}

\usepackage{textcmds} 
\usepackage{amsmath}
\usepackage{amsxtra}
\usepackage{amscd}
\usepackage{amsthm}
\usepackage{amsfonts}
\usepackage{amssymb}
\usepackage{eucal}
\usepackage[all]{xy}
\usepackage{graphicx}
\usepackage{comment}
\usepackage{epsfig}
\usepackage{psfrag}
\usepackage{mathrsfs}
\usepackage{amscd}
\usepackage{rotating}
\usepackage{lscape}
\usepackage{amsbsy}
\usepackage{verbatim}
\usepackage{moreverb}
\usepackage{url}
\usepackage{extarrows}
\usepackage{setspace}

\usepackage{cancel}
\usepackage{bbm}

\usepackage[hidelinks]{hyperref}

\usepackage{tikz-cd}
\usepackage{tikz}
\usetikzlibrary{matrix,arrows,decorations.pathmorphing}

\newtheorem{mainthm}{Theorem}

\newenvironment{mainthmbis}[1]
  {\renewcommand{\themainthm}{\ref{#1}$'$}%
   \addtocounter{mainthm}{-1}%
   \begin{mainthm}}
  {\end{mainthm}}

\newtheorem{cor}[subsubsection]{Corollary}
\newtheorem{lem}[subsubsection]{Lemma}
\newtheorem{prop}[subsubsection]{Proposition}

\newtheorem{thm}[subsubsection]{Theorem}

\newtheorem{defn}[subsubsection]{Definition}

\newtheorem{warning}[subsubsection]{Warning}

\theoremstyle{remark}
\newtheorem{rem}[subsubsection]{Remark}
\newtheorem{example}[subsubsection]{Example}

\numberwithin{equation}{section}

\newcommand{\nc}{\newcommand}
\nc{\renc}{\renewcommand}
\nc{\ssec}{\subsection}
\nc{\sssec}{\subsubsection}
\nc{\on}{\operatorname}

\nc\ol{\overline}
\nc{\ul}{\underline}
\nc\wt{\widetilde}

\renc{\setminus}{\smallsetminus}

\renewcommand{\subset}{\subseteq}
\renewcommand{\supset}{\supseteq}

\DeclareMathOperator{\Vect}{Vect} 
 
\DeclareMathOperator{\Top}{Top} 
\DeclareMathOperator{\Spc}{Spc}
 
\DeclareMathOperator{\Hom}{Hom} 
\DeclareMathOperator{\Sym}{Sym} 
\DeclareMathOperator{\Maps}{Maps} 
\DeclareMathOperator{\fSet}{fSet} 

\DeclareMathOperator{\oblv}{oblv} 

\DeclareMathOperator*{\colim}{colim}
\DeclareMathOperator*{\coker}{coker}
\DeclareMathOperator*{\cone}{cone}
\DeclareMathOperator*{\utimes}{\times}

\DeclareMathOperator{\PreStk}{PreStk} 
\DeclareMathOperator{\Aff}{Aff} 
\DeclareMathOperator{\Spec}{Spec} 
\DeclareMathOperator{\pt}{pt} 
\DeclareMathOperator{\Dom}{Dom}  
\DeclareMathOperator{\Sect}{Sect} 

\DeclareMathOperator{\QSect}{QSect} 
\DeclareMathOperator{\GMaps}{GMaps} 
\DeclareMathOperator{\GSect}{GSect} 

\DeclareMathOperator{\Par}{Par} 
\DeclareMathOperator{\Nilp}{Nilp} 
\DeclareMathOperator{\Whit}{Whit} 
\DeclareMathOperator{\Bun}{Bun} 
\DeclareMathOperator{\Ran}{Ran} 
\DeclareMathOperator{\rk}{rk} 
\DeclareMathOperator{\LS}{LS} 
\DeclareMathOperator{\coeff}{coeff} 
\DeclareMathOperator{\St}{St} 
\DeclareMathOperator{\Op}{Op} 
\DeclareMathOperator{\GFlags}{GFlags} 

\nc{\Opbar}{\ol\Op}
\nc{\Opbargen}{\ol\Op^{\gen}}
\nc{\Opgen}{\Op^{\gen}}

\DeclareMathOperator{\Coh}{Coh} 
\DeclareMathOperator{\Perf}{Perf} 
\DeclareMathOperator{\QCoh}{QCoh} 
\DeclareMathOperator{\IndCoh}{IndCoh} 
\DeclareMathOperator{\ICoh}{IndCoh} 
\DeclareMathOperator{\Dmod}{D-mod} 

\newcommand{\BM}{\mathrm{BM}} 
\newcommand{\dR}{\mathrm{dR}} 
\newcommand{\op}{\mathrm{op}} 
\let\ss\relax
\DeclareMathOperator{\ss}{ss} 
\newcommand{\ext}{\mathrm{ext}} 
\newcommand{\inc}{\mathrm{inc}} 
\newcommand{\indep}{\mathrm{indep}} 
\newcommand{\triv}{\mathrm{triv}} 

\newcommand{\gen}{\mathrm{gen}} 
\newcommand{\ggen}{\mathrm{-gen}} 

\newcommand{\temp}{\mathrm{temp}} 

\nc{\lambdach}{{\check\lambda}}
\nc{\Lambdach}{{\check\Lambda}{}}
\nc{\much}{{\check\mu}}
\nc{\omegach}{{\check\omega}}
\nc{\nuch}{{\check\nu}}
\nc{\etach}{{\check\eta}}
\nc{\alphach}{{\check\alpha}}

\nc{\rhoch}{{\check\rho}}

\nc{\cLambda}{{\check\Lambda}}

\nc{\cla}{{\check\lambda}}
\nc{\cmu}{{\check\mu}}
\nc{\cnu}{{\check\nu}}
\nc{\ceta}{{\check\eta}}

\nc{\bbone}{\mathbbm{1}}
\nc{\bbA}{{\mathbb{A}}}
\nc{\bbB}{\mathbb{B}}
\nc{\bbC}{{\mathbb{C}}}
\nc{\bbD}{{\mathbb{D}}}
\nc{\bbE}{{\mathbb{E}}}
\nc{\bbF}{{\mathbb{F}}}
\nc{\bbG}{{\mathbb{G}}}
\nc{\bbH}{{\mathbb{H}}}
\nc{\bbI}{{\mathbb{I}}}
\nc{\bbJ}{{\mathbb{J}}}
\nc{\bbK}{{\mathbb{K}}}
\nc{\bbL}{{\mathbb{L}}}
\nc{\bbM}{{\mathbb{M}}}
\nc{\bbN}{{\mathbb{N}}}
\nc{\bbO}{{\mathbb{O}}}
\nc{\bbP}{{\mathbb{P}}}
\nc{\bbQ}{{\mathbb{Q}}}
\nc{\bbR}{{\mathbb{R}}}
\nc{\bbS}{{\mathbb{S}}}
\nc{\bbT}{{\mathbb{T}}}
\nc{\bbU}{{\mathbb{U}}}
\nc{\bbV}{{\mathbb{V}}}
\nc{\bbW}{{\mathbb{W}}}
\nc{\bbX}{{\mathbb{X}}}
\nc{\bbY}{{\mathbb{Y}}}
\nc{\bbZ}{{\mathbb{Z}}}

\nc{\CA}{{\mathcal{A}}}
\nc{\CB}{{\mathcal{B}}}

\nc{\CE}{{\mathcal{E}}}
\nc{\CF}{{\mathcal{F}}}
\nc{\CH}{{\mathcal{H}}}

\nc{\CL}{{\mathcal{L}}}
\nc{\CC}{{\mathcal{C}}}
\nc{\CG}{{\mathcal{G}}}
\nc{\CM}{{\mathcal{M}}}
\nc{\CN}{{\mathcal{N}}}
\nc{\CK}{{\mathcal{K}}}
\nc{\CO}{{\mathcal{O}}}
\nc{\CP}{{\mathcal{P}}}
\nc{\CQ}{{\mathcal{Q}}}
\nc{\CR}{{\mathcal{R}}}
\nc{\CS}{{\mathcal{S}}}
\nc{\CU}{{\mathcal{U}}}
\nc{\CV}{{\mathcal{V}}}
\nc{\CW}{{\mathcal{W}}}
\nc{\CX}{{\mathcal{X}}}
\nc{\CY}{{\mathcal{Y}}}
\nc{\CZ}{{\mathcal{Z}}}
\nc{\CI}{{\mathcal{I}}}
\nc{\I}{\CI}

\nc{\fa}{{\mathfrak{a}}}
\nc{\fb}{{\mathfrak{b}}}
\nc{\fc}{{\mathfrak{c}}}
\nc{\fd}{{\mathfrak{d}}}
\nc{\ff}{{\mathfrak{f}}}
\nc{\fg}{{\mathfrak{g}}}
\nc{\fgl}{{\mathfrak{gl}}}
\nc{\fh}{{\mathfrak{h}}}
\nc{\fj}{{\mathfrak{j}}}
\nc{\fl}{{\mathfrak{l}}}
\nc{\fm}{{\mathfrak{m}}}
\nc{\fn}{{\mathfrak{n}}}
\nc{\fu}{{\mathfrak{u}}}
\nc{\fp}{{\mathfrak{p}}}
\nc{\fr}{{\mathfrak{r}}}
\nc{\fs}{{\mathfrak{s}}}
\nc{\ft}{{\mathfrak{t}}}
\nc{\fz}{{\mathfrak{z}}}
\nc{\fsl}{{\mathfrak{sl}}}
\nc{\hsl}{{\widehat{\mathfrak{sl}}}}
\nc{\hgl}{{\widehat{\mathfrak{gl}}}}
\nc{\hg}{{\widehat{\mathfrak{g}}}}
\nc{\chg}{{\widehat{\mathfrak{g}}}{}^\vee}
\nc{\hn}{{\widehat{\mathfrak{n}}}}
\nc{\chn}{{\widehat{\mathfrak{n}}}{}^\vee}

\nc{\fA}{{\mathfrak{A}}}
\nc{\fB}{{\mathfrak{B}}}
\nc{\fD}{{\mathfrak{D}}}
\nc{\fE}{{\mathfrak{E}}}
\nc{\fF}{{\mathfrak{F}}}
\nc{\fG}{{\mathfrak{G}}}
\nc{\fK}{{\mathfrak{K}}}
\nc{\fL}{{\mathfrak{L}}}
\nc{\fM}{{\mathfrak{M}}}
\nc{\fN}{{\mathfrak{N}}}
\nc{\fP}{{\mathfrak{P}}}
\nc{\fR}{{\mathfrak{R}}}
\nc{\fS}{{\mathfrak{S}}}
\nc{\fU}{{\mathfrak{U}}}
\nc{\fV}{{\mathfrak{V}}}
\nc{\fX}{{\mathfrak{X}}}
\nc{\fY}{{\mathfrak{Y}}}
\nc{\fZ}{{\mathfrak{Z}}}

\nc{\bb}{{\mathbf{b}}}
\nc{\bc}{{\mathbf{c}}}
\nc{\bd}{{\mathbf{d}}}
\nc{\bbf}{{\mathbf{f}}}
\nc{\be}{{\mathbf{e}}}
\nc{\bg}{{\mathbf{g}}}
\nc{\bi}{{\mathbf{i}}}
\nc{\bj}{{\mathbf{j}}}
\nc{\bn}{{\mathbf{n}}}
\nc{\bo}{{\mathbf{o}}}
\nc{\bp}{{\mathbf{p}}}
\nc{\bq}{{\mathbf{q}}}
\nc{\bt}{{\mathbf{t}}}
\nc{\bu}{{\mathbf{u}}}
\nc{\bv}{{\mathbf{v}}}
\nc{\bx}{{\mathbf{x}}}
\nc{\bs}{{\mathbf{s}}}
\nc{\by}{{\mathbf{y}}}
\nc{\bw}{{\mathbf{w}}}
\nc{\bA}{{\mathbf{A}}}
\nc{\bK}{{\mathbf{K}}}
\nc{\bB}{{\mathbf{B}}}
\nc{\bC}{{\mathbf{C}}}
\nc{\bG}{{\mathbf{G}}}
\nc{\bD}{{\mathbf{D}}}
\nc{\bH}{{\mathbf{H}}}
\nc{\bM}{{\mathbf{M}}}
\nc{\bN}{{\mathbf{N}}}
\nc{\bO}{{\mathbf{O}}}
\nc{\bT}{{\mathbf{T}}}
\nc{\bV}{{\mathbf{V}}}
\nc{\bW}{{\mathbf{W}}}
\nc{\bX}{{\mathbf{X}}}
\nc{\bZ}{{\mathbf{Z}}}
\nc{\bS}{{\mathbf{S}}}

\nc{\sA}{{\mathsf{A}}}
\nc{\sB}{{\mathsf{B}}}
\nc{\sC}{{\mathsf{C}}}
\nc{\sD}{{\mathsf{D}}}
\nc{\sF}{{\mathsf{F}}}
\nc{\sG}{{\mathsf{G}}}
\nc{\sK}{{\mathsf{K}}}
\nc{\sM}{{\mathsf{M}}}
\nc{\sO}{{\mathsf{O}}}
\nc{\sW}{{\mathsf{W}}}
\nc{\sQ}{{\mathsf{Q}}}
\nc{\sP}{{\mathsf{P}}}
\nc{\sR}{{\mathsf{R}}}
\nc{\sS}{{\mathsf{S}}}
\nc{\sT}{{\mathsf{T}}}
\nc{\sV}{{\mathsf{V}}}
\nc{\sZ}{{\mathsf{Z}}}

\nc{\sfp}{{\mathsf{p}}}
\nc{\sll}{{\mathsf{l}}}
\nc{\sr}{{\mathsf{r}}}
\nc{\bk}{{\mathsf{k}}}
\nc{\sg}{{\mathsf{g}}}

\nc{\sff}{{\mathsf{f}}}
\nc{\sfb}{{\mathsf{b}}}
\nc{\sfc}{{\mathsf{c}}}
\nc{\sd}{{\mathsf{d}}}
\nc{\se}{{\mathsf{e}}}

\nc{\one}{{\mathbf{1}}}
\nc{\two}{{\mathbf{t}}}

\nc{\cT}{{\check{T}}}
\nc{\cG}{{\check{G}}}
\nc{\cM}{{\check{M}}}
\nc{\cB}{{\check{B}}}
\nc{\cP}{{\check{P}}}

\nc{\ct}{{\check{\mathfrak t}}}
\nc{\cg}{{\check{\fg}}}
\nc{\cb}{{\check{\fb}}}
\nc{\cn}{{\check{\fn}}}
\nc{\cp}{{\check{\fp}}}
\nc{\cm}{{\check{\fm}}}

\nc{\inftyGrpd}{{\mathsf{Grpd}_\infty}}

\renc{\O}{\CO}
\nc{\Y}{\CY}

\nc{\restr}[2]{\left. #1 \right |_{#2}}
\nc{\uprestr}[2]{\left. #1 \right |^{#2}}

\nc{\hto}{\hookrightarrow}
\nc{\xto}{\xrightarrow}

\renc{\sec}{\section}

\nc{\lr}{\xymatrix{ \ar@<-0.4ex>[r] \ar@<.5ex>[l]  & } }
\nc{\rr}{\xymatrix{ \ar@<-0.2ex>[r] \ar@<.7ex>[r]  & } }
\nc{\rrr}{\xymatrix{ \ar@<.0ex>[r] \ar@<.7ex>[r] \ar@<-0.7ex>[r] & } }

\renc{\cong}{\simeq}

\nc{\virg}[1]{``#1"}

\renc{\bold}[1]{\boldsymbol{#1}}

\nc{\bigt}[1]{\big( #1 \big) }
\nc{\Bigt}[1]{\Big( #1 \Big) }

\nc{\footcite}{\footnote}

\nc{\GA}{{G(\bbA)}}
\nc{\GK}{{G(\bbK)}}
\nc{\GO}{{G(\bbO)}}
\nc{\NO}{{N(\bbO)}}

\nc{\NA}{N(\AA)}
\nc{\VA}{V(\AA)}

\nc{\squigto}{\rightsquigarrow}
\nc{\longto}{\longrightarrow}
\nc{\tto}{\twoheadrightarrow}

\nc{\Gch}{{\check{G}}}
\nc{\Pch}{{\check{P}}}
\nc{\Mch}{{\check{M}}}
\nc{\Qch}{{\check{Q}}}

\nc{\LL}{\mathbb{L}}

\renc{\hat}{\widehat}

\nc{\heart}{\heartsuit}
\nc{\kk}{\mathbbm{k}} 

\nc{\HHom}{\CH{om}} 

\nc{\gch}{\mathfrak{\check{g}}}

\nc{\LSGch}{{\LS_\Gch}}
\nc{\LSPch}{{\LS_\Pch}}
\nc{\LSMch}{{\LS_\Mch}}

\nc{\Hsx}[2]{\H_{{#1} \leftarrow {#2}}}
\nc{\Hdx}[2]{\H_{{#1} \to {#2}}}
\nc{\Hcorr}[3]{ \H_{{#1} \leftarrow {#2} \to {#3}} }
\nc{\Hopcorr}[3]{ \H_{{#1} \to {#2} \leftto {#3}} }

\nc{\ICohsx}[2]{\ICohW_{{#1} \leftarrow {#2}}}
\nc{\ICohdx}[2]{\ICohW_{{#1} \to {#2}}}
\nc{\ICohcorr}[3]{ \ICohW_{{#1} \leftarrow {#2} \to {#3}} }
\nc{\ICohopcorr}[3]{ \ICohW_{{#1} \to {#2} \leftto {#3}} }

\nc{\QCohsx}[2]{\QCohW_{{#1} \leftarrow {#2}}}
\nc{\QCohdx}[2]{\QCohW_{{#1} \to {#2}}}
\nc{\QCohcorr}[3]{ \QCohW_{{#1} \leftarrow {#2} \to {#3}} }
\nc{\QCohopcorr}[3]{ \QCohW_{{#1} \to {#2} \leftto {#3}} }

\nc{\leftto}{\leftarrow}
\nc{\lto}{\leftto}
\nc{\xlto}[1]{\xleftarrow{#1}}

\nc{\ltemp}{{}^\temp}
\nc{\lcusp}{{}^\cusp}
\nc{\TwCorr}{\mathsf{TwCorr}}

\nc{\aft}{{\mathit{aft}}}

\nc{\g}{\mathfrak{g}}

\nc{\LSG}{\LS_G}
\nc{\LSM}{\LS_M}
\nc{\LSP}{\LS_P}
\nc{\LST}{\LS_T}
\nc{\LSB}{\LS_B}
\nc{\Gm}{{\GG_m}}

\nc{\olBun}{\ol{\Bun}}
\nc{\cl}{\on{cl}}
\nc{\starext}{{*\on{-ext}}}

\nc{\RHom}{\on{RHom}}

\nc{\mon}{{\Gm\on{-mon}}}

\nc{\TBunG}{\sT_{\Bun_G}}
\nc{\lperp}{{}^\perp}

\nc{\Levi}{\on{Levi}}

\nc{\pscpt}{{\on{ps-cpt}}}

\nc{\Hren}{H^{{\BM}}_*}

\nc{\sigmatriv}{\sigma_{\triv}}

\nc{\QF}{\mathit{QF}}

\nc{\mRan}{\mathrm{Ran}_\ast}

\nc{\CHom}{\CH{om}}

\begin{document}

\title{Contractibility of the space of generic opers for classical groups}

\author{Dario Beraldo \and David Kazhdan \and Tomer M. Schlank}

\begin{abstract}
Let $G$ be a reductive group and $X$ a smooth projective curve.
We prove that, for $G$ classical and $\sigma$ an arbitrary $G$-local system on $X$, the space $\Opbargen_{G,\sigma}$ of generic extended oper structures on $\sigma$ is homologically contractible.
This contractibility result is crucial for the proof of the geometric Langlands conjecture.
\end{abstract}

\maketitle

\sec{Introduction}

In this paper, we prove two important results on the spectral side of the geometric Langlands conjecture: the first one (Theorem \ref{mainthm:contract-fiber}) is a contractibility statement for the space of generic Griffiths transverse flags on a $G$-local system, the second one (Theorem \ref{mainthm:St-Opers}) relates the space of generic oper structures to Deligne-Lusztig duality.

Our techniques are quite general and likely to be useful in other contexts. For instance, we prove the $\O$-contractibility of $\fL := \GMaps(X,\bbP^n)$, the space of generic maps from a curve to $\bbP^n$.
We also prove a surprising fact (Theorem \ref{mainthm:Dmod-in-ICoh}): the forgetful functor $\Dmod(\fL) \to \QCoh(\fL)$ is fully faithful.
These facts yield a new proof the homological contractibility of $\fL$.\footnote{The latter result was first proven in \cite{Barlev}, relying heavily on \cite{contract}. Our method is more direct.}
In the final part of the paper, we exploit our methods to prove a strengthening of Tsen's theorem (Theorem \ref{mainthm:Tsen}).

\ssec{The first main theorem}

\sssec{}

Let $k$ be  an algebraically closed field of characteristic $0$ (fixed for the remainder of the paper), $G$ a reductive $k$-group and $X$ a smooth projective $k$-curve. We denote by $\LSG:=\LSG(X)$ the moduli stack of de Rham $G$-local systems on $X$. A typical $k$-point of $\LSG$ will be denoted by $\sigma = (E,\nabla)$, where $E$ stands for the underlying $G$-bundle and $\nabla$ for the flat connection.

\sssec{}

Let us state our main result, deferring the definition of the terms involved to later paragraphs.

\begin{mainthm} \label{mainthm:contract-fiber}
Let $G$ be a classical reductive group, that is, a reductive group of Dynkin type A, B, C or D. For any $G$-local system $\sigma$ on $X$, the space $\Opbargen_{G,\sigma}$ of \emph{generic extended oper structures} on $\sigma$ is homologically contractible.
\end{mainthm}

\begin{rem}
For us, the expressions \virg{generic extended oper structures} and \virg{generic Griffiths transverse flags} are synonyms. We will mostly use the former expression.
The reason for the adjective \virg{extended} will be clarified in the sequel: it has been chosen in order to match with the \emph{extended Whittaker category} appearing on the automorphic side of the geometric Langlands conjecture. In fact, we will also consider the space $\Opgen_{G,\sigma}$ of (\emph{non-extended}, or \emph{genuine}) generic opers structures on $\sigma$. The latter space is an essential player in the geometric Laglands program and will be the subject of our second main result.
\end{rem}

\begin{rem}
We expect Theorem \ref{mainthm:contract-fiber} to hold for the exceptional groups as well. In this direction, D. Arinkin proved that $\Opbargen_{G,\sigma}$ is nonempty for any $G$, see \cite{Arinkin}.
\end{rem}

\sssec{}

Let us proceed with the definition of $\Opbargen_{G,\sigma}$. Fixing a Borel subgroup $B \subset G$ once and for all, denote by $N \subset B$ the unipotent radical and by $\fn \subseteq {\mathfrak b} \subseteq \g$ the corresponding inclusions of Lie algebras. 
Consider the vector space
$$
\g^{\geq -1} :=
\{x \in \g \, | \,
[x, \fn] \subset \fb
\}.
$$
Clearly, $\g^{\geq -1}$ is $B$-invariant and contains $\fb$ as a sub-representation. 
Moreover, the quotient $\g^{\geq -1}/\fb$ splits $B$-equivariantly as the direct sum of the negative simple root spaces:
$$
\g^{\geq -1}/\fb
\simeq
\bigoplus_{\alpha \in \sS} \g^{(-\alpha)}.
$$

\sssec{}

Let $\sigma = (E,\nabla)$ be a $G$-local system on a smooth (not necessarily proper) curve $X'$. Given a $B$-reduction $E_B$ of $E$, we can canonically obtain a global section $\nabla_{|E_B}$ of the vector bundle $(\g/\fb)_{E_B} \otimes \Omega_{X'}$. This section measures the discrepancy between the $B$-reduction and the connection: for instance, $\nabla_{|E_B} = 0$ iff $E_B$ is preserved by the connection, that is, iff $E_B$ is a $B$-reduction of $\sigma$ as a $G$-local system. 

Any $B$-invariant subspace contained between $\fb$ and $\fg$ gives rise to a sub-bundle of $(\g/\fb)_{E_B} \otimes \Omega_{X'}$. In our case of interest, we consider the subbundle
$$
(\g^{\geq -1}/\fb)_{E_B} \otimes \Omega_{X'}
\subseteq
(\g/\fb)_{E_B} \otimes \Omega_{X'}
$$
and give the following important definition.

\begin{defn}
Let $\sigma=(E,\nabla)$ be a $G$-local system on a smooth (not necessarily proper) curve $X'$.
We say that a $B$-reduction $E_B$ of $E$ is an \emph{extended oper structure} if the section $\nabla_{|E_B}$ belongs to the subbundle
$(\g^{\geq -1}/\fb)_{E_B} \otimes \Omega_{X'}$.
\end{defn}

\sssec{}

Let $X'$ be a smooth curve and $\sigma=(E,\nabla)$ a $G$-local system on $X'$.
Now we are ready to give the first (informal) definition of the space $\Opbargen_{G,\sigma}$. Consider the space of pairs $(U,E_{B,U})$, where:
\begin{itemize}
\item
$U$ is a nonempty open subset of $X'$;
\item
$E_{B,U}$ is a $B$-reduction of $\restr E U$, which satisfies the extended oper structure condition.
\end{itemize}
There is an obvious equivalence relation, which identifies $(U_1,E_{B,U_1}) \approx (U_2,E_{B,U_2})$ iff these two pieces of data coincide on a nonempty open subset of $U_1 \cap U_2$.
We let $\Opbargen_{G,\sigma}$ be the space that classifies equivalence classes of pairs $(U,E_{B,U})$ as above.

\sssec{}

A proper definition of $\Opbargen_{G, \sigma}$ will be provided in the main body of the paper.
It turns out that $\Opbargen_{G,\sigma}$, when defined rigorously as a functor of points, is represented by a \emph{pseudo-scheme}: a (not necessarily filtered) colimit of schemes of finite type along proper transition maps.

Now, any pseudo-scheme $\Y$ admits a well-defined differential graded (DG) category of D-modules $\Dmod(\Y)$, which furthermore comes with a canonical object $\omega_\Y$, the dualizing sheaf. We say that $\Y$ is \emph{homologically contractible} if the DG vector space 
$$
H^*(\Y):= \CHom_{\Dmod(\Y)}(\omega_\Y,\omega_\Y)
$$
is quasi-isomorphic to the ground field $k$. In Section 
\ref{ssec:prelim-contract}, we will briefly recall why this definition is reasonable.
As stated above, the main goal of this paper is to prove the homological contractibility, in the above sense, of $\Opbargen_{G,\sigma}$.

\ssec{Motivation}

Let us now explain our interest in such a contractibility and the role it plays in the geometric Langlands program. As a preliminary step, we need to slightly reformulate Theorem \ref{mainthm:contract-fiber}.

\sssec{}

Denote by $\Opbargen_G$ the prestack parametrizing triples $(E, \nabla; E_B)$ where:
\begin{itemize}
\item
$(E,\nabla)$ is a $G$-local system on the whole $X$;
\item
$E_B$ is a \emph{generic} $B$-reduction of $E$, which satifies the extended oper condition (where defined).
\end{itemize}
Let $\ol\pi: \Opbargen_G \to \LSG$ denote the forgetful map: by construction, its fiber over a $k$-point $\sigma$ of $\LSG$ is precisely the space $\Opbargen_{G,\sigma}$ considered above.

A simple argument will show that the following statement is equivalent to Theorem \ref{mainthm:contract-fiber}:

\begin{mainthmbis}{mainthm:contract-fiber} \label{mainthm:contract}
Let $G$ be a classical group. Then the functor
\begin{equation} \label{eqn-fully-faith-ext-opers}
\ol\pi^!:
\Dmod(\LSG) \longto
\Dmod(\Opbargen_G)
\end{equation}
is fully faithful. In other words, the natural arrow
\begin{equation} \label{eqn:nat-map-Opbar}
\ol\pi_!(\omega_{\Opbargen_G}) \longto \omega_{\LSG}
\end{equation}
is an isomorphism in $\Dmod(\LSG)$.
\end{mainthmbis}

\begin{rem}
The two statements in the above theorem are equivalent thanks to the pseudo-properness\footnote{See Section \ref{sssec:pseudo-proper map} for the definition.} of $\ol\pi$. Indeed, in this case, the left adjoint functor $\ol\pi_!$ is well-defined and satisfies the projection formula:
$$
\ol\pi_! \circ \ol\pi^!(-)
\simeq
\ol\pi_!(\omega_{\Opbargen_G}) \stackrel ! \otimes -.
$$
\end{rem}

\sssec{}

According to \cite{AG1}, the geometric Langlands conjecture amounts to an equivalence of DG categories
\begin{equation} \label{eqn:GL}
\ICoh_{\Nilp}(\LSG)
\stackrel ? \simeq
\Dmod(\Bun_{\Gch}),
\end{equation}
where $\Gch$ is the Langlands-dual of $G$ and $\ICoh_{\Nilp}$ is a certain enlargement of the DG category $\QCoh(\LSG)$. The strategy envisioned by D. Gaitsgory (see \cite{Outline}) hinges on having a fully faithful embedding of $\Dmod(\Bun_{\Gch})$ into a larger DG category $\Whit(\Gch,\ext)$, called the \emph{extended Whittaker category}.

\sssec{}

While the DG category $\Whit(\Gch,\ext)$ and the relevant functor
\begin{equation} \label{coeff-ext}
\coeff_{\Gch, \ext}:
\Dmod(\Bun_{\Gch})
\longto
\Whit(\Gch,\ext)
\end{equation}
have been defined (see \cite{Outline} and \cite{ext-whit}), it seems very difficult\footnote{Exceptions are the cases of $G = GL_n$ and $G=SL_n$, which were dealt with in \cite{ext-whit}.} to prove that $\coeff_{\Gch, \ext}$ is fully faithful. 
For this reason, we propose to \virg{vary the quantum parameter} from $0$ to $\infty$ and thus consider the analogous problem on the spectral side.
More precisely, the geometric Laglands equivalence \eqref{eqn:GL} should be the value at $0$ of a family of equivalences parametrized by $\bbP^1$, see \cite{quantum, Zhao}. 
The value at $\infty$ is the conjectural equivalence
\begin{equation} \label{eqn:GL-at-infty}
\Dmod(\Bun_{G})
\stackrel ? \simeq
\ICoh_{\Nilp}(\LSGch).
\end{equation}
It is natural to expect that the functor \eqref{coeff-ext} deforms accordingly: at $\infty$, it ought to yield a version of the fully faithfulness statement close to our Theorem \ref{mainthm:contract}.
Indeed, the similarities between extended/genuine Whittaker coefficients and extended/genuine opers are evident.

\sssec{}

Of course, one might object that there is a discrepancy between the fully faithfulness of \eqref{coeff-ext} and that of \eqref{eqn-fully-faith-ext-opers}: the latter should involve $\IndCoh_{\Nilp}(\LSG)$ and not $\Dmod(\LSG)$.
We can respond to this objection in two ways:
\begin{itemize}

\item
there is a way to modify Gaitsgory's strategy so that the homological contractibility of Theorem \ref{mainthm:contract} suffices, see \cite{ramanujan} for some details;

\item

Alternatively, it is possible to adapt the arguments of the present paper to prove the  $\O$-contractibility of the fibers of $\ol \pi$. We partially do this in Section \ref{sec:contract-LE}, where we prove that $\fL := \GMaps(X,\bbP ^n)$ is both $\O$-contractible and homologically contractible.
\end{itemize}

\begin{rem}
We do not pursue this $\O$-contractibility in full generality for two reasons. First, it would require using some derived algebraic geometry (while we wish to keep this paper readable to a larger audience); secondly, in view of the first item above, it might not be necessary.
\end{rem}

\ssec{Genuine opers}

So far, we have discussed the notion of extended opers. Now let us recall the related notion of (non-extended) opers, which actually arose first in the literature\footnote{see, for instance, \cite{BD-quant, Frenkel} and references therein} and which plays an essential role in our story as well.

\sssec{}

Let $\Opgen_G$ the prestack parametrizing triples $(E, \nabla; E_B)$ where:
\begin{itemize}
\item
$(E,\nabla)$ is a $G$-local system on the whole $X$;
\item
$E_B$ is a \emph{generic} $B$-reduction of $E$, which satifies the extended oper condition (where defined);
\item
we impose yet another condition. Observe that the vector bundle $(\g^{\geq -1}/\fb)_{E_B}$ splits as a direct sum of line bundles
$$
(\g^{\geq -1}/\fb)_{E_B}
 \simeq 
 \bigoplus_{\alpha \in \sS}(E_B)^{(-\alpha)},
 $$
one for each negative simple root; accordingly, the section $\nabla_{|E_B} \simeq \oplus_{\alpha}\nabla_{|E_B}^{-\alpha}$ is a direct sum of sections of these line bundles. We require all these sections to be generically nonzero.

\end{itemize}
\sssec{}

The natural inclusion $j: \Opgen_G \to \Opbargen_G$ is an open embedding, and it is an isomorphism over the locus of irreducible $G$-local systems. Indeed, the vanishing of one of the $\nabla_{|E_B}^{-\alpha}$ would yield a reduction of the $G$-local system in question to the maximal parabolic determined by $\alpha$, thereby contradicting irreducibility.
Let $\pi: \Opgen_G \to \LSG$ denote the natural forgetful map (which, by construction, equals $ \ol\pi \circ j$). 
Our next theorem concerns the D-module $\pi_{*,\dR}(\omega_{\Opgen_G})$. To state our result, we need a small digression on the Deligne-Lusztig duality functor.

\ssec{The second main theorem}

\sssec{}

The study of Deligne-Lusztig duality in the setting of the geometric Langlands program (see \cite{DL-aut} for the automorphic side and \cite{DL-spec} for the spectral side) has been particularly useful in generating provable conjectures, see for instance \cite{ramanujan} and \cite{aut-gluing}.

\sssec{}

Without going into details, let us recall that the study on the spectral side led the first author to the introduction of the \emph{Steinberg} D-module 
$$
\St_G \in \Dmod(\LSG).
$$ 
This is the D-module on $\LSG$ defined by the formula
\begin{equation} \label{eqn:Steinberg-defn}
\St_G
:=
\coker
\Bigt{
\colim_{P \in \Par'}
(\fp_P)_!(\omega_{\LSP})
\to
\omega_{\LSG}
},
\end{equation}
where $\Par' = \{P \, |\, B \subseteq P \subsetneq G \}$ denotes the poset of standard parabolics with $G$ excluded, $\fp_P: \LSP \to \LSG$ is the (proper) induction map, and the arrows forming the colimit are the obvious ones coming from adjunction.

\sssec{}

Here are three features of $\St_G$ that make it special.
\begin{itemize}
\item
The functor of action by $\St_G$ on $\QCoh(\LS_G)$ yields a fully faithful embedding
$$
\Perf(\LSG)
\hto
\QCoh(\LSG),
$$
which is different from the tautological embedding.

\item

The $k$-fibers of $\St_G$ are either $1$-dimensional (in some cohomological degree) or zero, according to whether the $G$-local system in question is semisimple or not.
For instance, the fiber $\restr{\St_G}{\sigma_\triv}$ at the trivial $G$-local system equals $k[|\sR| + \rk_G]$, where $\sR$ is the set of roots of $G$ and $\rk_G$ its semisimple rank.
At the other extreme,  $\restr{\St_G}{\sigma} \simeq k$ whenever $\sigma$ is irreducible.

\item

The principal monoidal ideal $\Dmod(\LSG)^{\St}$ generated by $\St_G$ admits a nice description, related to the fact that $\St_G$ can be thought as the dualizing sheaf of the ill-defined substack $\LSG^{\ss} \subseteq \LSG$ of \emph{semi-simple $G$-local systems}.

\end{itemize}
For reference, the above properties are \cite[Theorem E]{DL-spec}, \cite[Theorem D']{DL-spec} and \cite[Theorem C']{DL-spec}, respectively.

\sssec{}

Our second main result states that $\St_G$ admits a very natural description in terms of generic (genuine) opers, as follows. Recall from above the prestack $\Opgen_G$ of generic opers and the structure map $\pi: \Opgen_G \to \LSG$.

\begin{mainthm} \label{mainthm:St-Opers}
For $G$ a classical group, there is a canonical isomorphism 
$$
\pi_{*,\dR}(\omega_{\Opgen_G}) \simeq \St_G
$$
in $\Dmod(\LSG)$.
\end{mainthm}

\begin{rem}
Since $\pi$ is not pseudo-proper, the functors $\pi_{*,\dR}$ and $\pi_!$ are different. We do not know of a nice description of $\pi_!(\omega_{\Opgen_G})$.
\end{rem}

\begin{rem}
Part of the content of the above theorem is the construction of a natural arrow 
$$
\pi_{*,\dR}(\omega_{\Opgen_G}) \to \St_G.
$$
This issue does not occur in the setting of extended opers: in that case, the map \eqref{eqn:nat-map-Opbar} arises canonically from the counit of an adjunction.
\end{rem}

\sssec{}

Now let $\sigma$ be a $G$-local system defined on $X$, that is, a $k$-point of $\LSG$. Denote by $\Opgen_{G,\sigma}$ be the space of generic oper structures on $\sigma$.
Pulling back along $\sigma: \Spec(k) \to \LSG$ and applying base-change, we obtain information on the Borel-Moore homology of $\Opgen_{G,\sigma}$.

\begin{cor}
Retaining the same notation as above, for classical $G$, we have a canonical isomorphism
$$
\Hren(\Opgen_{G,\sigma}) \simeq \restr{\St_G}{\sigma}.
$$
Thus, $\Hren(\Opgen_{G,\sigma})$ is zero when $\sigma$ is not semi-simple and $1$-dimensional (placed in some cohomological degree) otherwise.
\end{cor}

\begin{example}
When $\sigma$ irreducible, the open embedding $\Opgen_{G,\sigma} \hto \Opbargen_{G,\sigma}$ is an isomorphism. Since $\Opbargen_{G,\sigma}$ is proper, we obtain that
$$
H_*(\Opgen_{G,\sigma})
\simeq 
\Hren(\Opgen_{G,\sigma}).
$$
On the other side, we already mentioned that $\restr {\St_G} \sigma \simeq k$. Thus we have recovered a particular case of Theorem \ref{mainthm:contract-fiber}. This is hardly a surprise: our strategy consists of deducing Theorem \ref{mainthm:St-Opers} from Theorem \ref{mainthm:contract}.
\end{example}

\begin{example}
Three more examples for the group $SL_2$ are illustrated in \cite{ramanujan}.
\end{example}

\sssec{}

As an extra consequence, we can partially answer a question left open in \cite{Arinkin} and \cite{Frenkel-Zhu}: the question is whether any $G$-local system on $X$ admits a generic oper structure.

\begin{cor}
Let $G$ be a classical reductive group and $\sigma \in \LSG(k)$ semi-simple. Then $\Opgen_{G,\sigma}$ is not empty. 
\end{cor}

\ssec{Outline of the proofs}

For the reader's convenience, let us give a quick account of the main ideas in the proofs of Theorem \ref{mainthm:contract-fiber} and Theorem \ref{mainthm:St-Opers}.

\sssec{}

First, Theorem \ref{mainthm:contract-fiber} and Theorem \ref{mainthm:contract} are easily seen to be equivalent. This follows from the fact that relative homological contractibility of a pseudo-proper map can be checked on field-valued fibers, see \cite[Lemma 6.1.7]{AG2} or \cite[Section 4.2]{epiga}.

\sssec{}

Next, we will deduce Theorem \ref{mainthm:St-Opers} from Theorem \ref{mainthm:contract}, the latter applied to $G$ and its standard Levi subgroups. Our argument is largely formal\footnote{In particular, the validity of Theorem \ref{mainthm:contract} for the exceptional groups would immediately imply that of Theorem \ref{mainthm:St-Opers} for the same groups.} and goes roughly as follows.
By definition, the difference between $\St_G$ and $\omega_{\LSG}$ is controlled by the proper standard parabolics:
$$
\ker
(\omega_{\LSG} \to \St_G)
\simeq
\colim_{P \in \Par'} \,
(\fp_P)_!(\omega_{\LSP}).
$$
Observe now that the difference between $\omega_{\Opbargen_G}$ and $j_*(\omega_{\Opgen_G})$ is controlled by $\Par'$, too: indeed, such a difference comes from the vanishing pattern of the sections $\nabla^\alpha_{| E_B}$, which in turn is determined by the non-trivial subsets of the Dynkin diagram.
Translating this into a formula (with notations explained in Section \ref{sec:Steinberg}), we will obtain an isomorphism
$$
\ker
\Bigt{
\omega_{\Opbargen_G} \to j_{*,\dR}(\omega_{\Opgen_G})
}
\simeq
\colim_{P \in \Par'} \,
(i_P)_! \bigt{ \omega_{\Opbargen_{G,P}} }.
$$
Theorem \ref{mainthm:contract} then allows to compare these two displayed expressions and thereby produce an isomorphism $\pi_{*,\dR}(\omega_{\Opgen_G}) \simeq \St_G$ as desired.

\sssec{}

It remains to prove Theorem \ref{mainthm:contract-fiber}. At this point, we assume $G$ to be a classical group and proceed for each type separately. In this introduction, we only carry out the example of $G=SO_7$; the other cases are similar and will be discussed in Section \ref{sec:classical}.

\sssec{}

An $SO_7$-local system $\sigma$ is a quadruple $(E,\beta, \nabla, \tau)$ where:
\begin{itemize}
\item
 $E$ is vector bundle $E$ of rank $7$;
\item
$\beta: Sym^2(E) \to \O_X$ is a non-degenerate symmetric bilinear form on $E$;
\item
$\nabla: E \to E \otimes \Omega_X$ is a connection on $E$, which is compatible with $\beta$: for any $\mu \in \bbT_X$ and $e_1, e_2 \in E$, we have
\begin{equation} \label{eqn:beta-vs-nabla}
\mu(\beta(e_1, e_2))
\simeq
\beta (\nabla_\mu(e_1), e_2) 
+ 
\beta (e_1, \nabla_\mu(e_2)); 
\end{equation}
\item
$\tau$ is a horizontal trivialization of the rank-$1$ local system $(\det(E), \nabla_{\det(E)})$.
\end{itemize}

\sssec{}

Note that $\Opbargen_{G,\sigma}$ depends only on the restriction of $\sigma$ at the generic point on $X$. So we can assume that both $E$ and $\Omega_X$ are trivialized (the latter by a vector field $\nu$): in this case, $\nabla$ is determined by the differential operator $\nabla_\nu: E \to E$.

Unraveling the definition, we see that a generic extended oper structure on $\sigma$ consists of:
\begin{itemize}
\item
a generically defined $SO_7$-flag on $E$, that is, a full isotropic flag of vector sub-bundles
$$
0 \subset E_1 \subset E_2 \subset E_3 \subset \restr E U = (\O_U)^{\oplus 7}, 
\hspace{.5cm}
\beta(E_3, E_3 )= 0,
$$
defined on a nonempty open subset $U \subseteq X$;

\item
the \emph{Griffiths transversality condition}, which in this case amounts to:
$$
\nabla_\nu(E_1) \subset E_2,
\hspace{.4cm}
\nabla_\nu(E_2) \subset E_3.
$$
\end{itemize}

\sssec{}

Now we come to the main idea. We consider an auxiliary pseudo-scheme $\Y_{G,\sigma}$, equipped with a map $\xi: \Opbargen_{G,\sigma} \to \Y_{G,\sigma}$. Elements of $\Y_{G,\sigma}$ consist of:
\begin{itemize}
\item
a generic line sub-bundle
$$
L \subset \restr E U \simeq \O_U^{\oplus 7}
$$
defined on a nonempty open subset $U \subseteq X$;

\item
the following five conditions:
\begin{eqnarray}  \label{eqn:nablas-SO(7)}
& & 
\beta(L,L ) = 0 
\\
\nonumber
& & 
\beta(L,\nabla_\nu(L)) = 0   
\\
\nonumber
& & 
\beta(\nabla_\nu(L),\nabla_\nu(L)) = 0
\\
\nonumber
& & 
\beta(\nabla_\nu(L),\nabla_\nu^2 (L)) = 0
\\
\nonumber
& & 
\beta(\nabla_\nu^2(L),\nabla_\nu^2(L)) = 0.
\end{eqnarray}

\end{itemize}

\sssec{}

These conditions evidently mean that the vector sub-bundle generated by 
$$
L, \nabla_\nu(L), \nabla_\nu^2(L)
$$
is $\beta$-isotropic. In particular, we have a well-defined map 
$$
\xi: \Opbargen_{G,\sigma} \to \Y_{G,\sigma}
$$
$$
(E_1 \subset E_2 \subset E_3)
\mapsto E_1.
$$
This map is pseudo-proper, hence the homological contractibility of $\Opbargen_{G,\sigma}$ can be split in two parts: 
\begin{itemize}
\item
the natural arrow
$
\xi_!(\omega_{\Opbargen_{G,\sigma}})
\to \omega_{\Y_{G,\sigma}}
$
is an isomorphism;

\item

the space $\Y_{G,\sigma}$ is homologically contractible.

\end{itemize}

\sssec{}

The first item is dealt with fiberwise. 
We will show that the nature of the fibers of $\xi$ depends on the generic rank of $\langle L, \nabla_\nu(L), \nabla_\nu^2(L) \rangle$: when the latter equals $i$ (with $1 \leq i \leq 3$), the fiber is isomorphic to $\Opbargen_{SO
_{2i-1}, \sigma'}$ for some $\sigma' \in \LS_{SO_{2i-1}}$. 
In particular, the fibers of $\xi$ are all homologically contractible by induction.

\sssec{}

It remains to prove that $\Y_{G,\sigma}$ is itself homologically contractible. To this end, we observe that the second and fourth equations in \eqref{eqn:nablas-SO(7)} are redundant: the second is implied by the first, the fourth by the third. Hence we have
$$
\Y_{G,\sigma}
\simeq
\left\{
L \in \GMaps(X,\bbP^6) 
\,
\big|
\,
\beta(L,L)
=
\beta(\nabla_\nu(L),\nabla_\nu(L))
=
\beta(\nabla_\nu^2(L),\nabla_\nu^2(L)) = 0
\right\},
$$ 
where $ \GMaps(X,\bbP^6) $ is the space of generic maps $X \dasharrow \bbP^6$. Thus, $\Y_{G,\sigma}$ is cut out in $\GMaps(X,\bbP^6) $ by \emph{three quadratic} equations. 
Then the homological contractibility of $\Y_{G,\sigma}$ is an instance of a general Lefschetz-type result, which is valid whenever the sum of the degrees of the equations does not the exceed the dimension of the ambient projective space: in our case $6 \geq 3 \cdot 2$.

\sssec{}

We warn the reader that the definition of $\Y_{G,\sigma}$ heavily depends on the Dynkin type of $G$. The relevant definitions for all four families appear in Section \ref{sec:classical}.

In general, for type $B_n$ we will have $n$ quadratic equations in $\GMaps(X,\bbP^{2n})$.
The situations for type $A_n$, $C_n$ and $D_n$ are even better:
\begin{itemize}
\item
in type $A_n$, we get no equations: in this case, $\Y_{G, \sigma}$ is simply $\GMaps(X,\bbP^n)$;

\item
in type $C_n$ and $D_n$, we get $n-1$ quadratic equations in $\GMaps(X,\bbP^{2n-1})$.
\end{itemize}

\begin{rem}

We failed to get a suitable set of equations for the exceptional groups. For example, for $G=G_2$ one can use the analysis in \cite{An09} to produce two quadratic equations and one cubic equation in $\GMaps(X,\bbP^6)$. Unfortunately, as $2+2+3 > 6$, our strategy cannot be applied.

\end{rem}

\sssec{}

The starting point for our Lefschetz-type estimates is the homologically contractibility of $\fL:= \GMaps(X,\bbP^n)$. We will give a new proof of this important fact. Our proof is stronger than the original one (see \cite{Barlev}) in that it establishes the $\O$-contractibility, as well as a surprising categorical statement:

\begin{mainthm} \label{mainthm:Dmod-in-ICoh}
The prestack $\fL$ is $\O$-contractible and the natural forgetful functor
$$
\oblv: \Dmod(\fL) 
\longto
\QCoh(\fL)
$$
is fully faithful. In particular, $H^*(\fL) \simeq H^*(\fL, \O) \simeq k$.
\end{mainthm}

\sssec{}

As mentioned in the beginning, our techniques are general and might turn out to be useful in other contexts: as a simple illustration of this, we show at the end of the paper that our methods yield a strengthening of Tsen's theorem.

\begin{mainthm} \label{mainthm:Tsen}
Let $\Y \subseteq \bbP^n_X$ be a closed subscheme defined by $m$ homogeneous equations of degrees $d_1, \ldots, d_m$. Consider the pseudo-scheme $\GSect(X, \Y)$ of generic sections of $\Y \to X$. If $\sum_i d_i \leq n$, then $\GSect(X, \Y)$ is homologically contractible.
\end{mainthm}

\ssec{Structure of the paper}
We begin in Section \ref{sec:prelim} by collecting many of the basic definitions and facts we need: some general algebraic geometry (pseudo-schemes and pseudo-proper maps, contractibility, the Ran space), as well as the notions of generic maps/sections and the definition of (extended) oper structures.
In Section \ref{sec:Steinberg}, we show how the Steinberg object arises from the space of generic genuine opers, thereby proving that Theorem \ref{mainthm:contract-fiber} implies Theorem \ref{mainthm:St-Opers}.
In Section \ref{sec:classical}, we focus of classical groups and start the proof of Theorem \ref{mainthm:contract-fiber}. Namely, we describe how $\Opbargen_{G,\sigma}$ looks like in each case and define an auxiliary space $\Y_{G,\sigma}$ with the same homology as $\Opbargen_{G,\sigma}$.
The tow sections afterwards are devoted to proving that $\Y_{G,\sigma}$ is homologically contractible: in Section \ref{sec:contract-LE}, we describe $\fL_E$ using the Ran space and prove Theorem \ref{mainthm:Dmod-in-ICoh}; in Section \ref{sec:Oka}, we use some t-structure estimates to show that $\Y_{G,\sigma}$ has the same homology as $\fL_E$.
Finally, in Section \ref{sec:Tsen}, we use the same t-structure estimates to prove Theorem \ref{mainthm:Tsen}.

\ssec{Acknowledgements}
We would like to thank Dennis Gaitsgory and
Lior Yanovsky for useful discussions, and the referee for many useful remarks and corrections. DB further thanks Ian Grojnowski, Massimo Pippi, Michael Singer, Bertrand To\"en and Gabriele Vezzosi for interesting conversations and support.

The research of DK has been partially supported by the ERC grant No 669655.
TMS is supported by ISF1588/18 and BSF 2018389.

\sec{Notation, background and preliminaries} \label{sec:prelim}

\ssec{Algebraic geometry and category theory}

\sssec{}

Throughout this paper, all schemes and prestacks (see below) we consider are defined over $k$, unless otherwise stated, and they are locally of finite type.
Accordingly, given two schemes $S$ and $T$, their product over $\Spec k$ will be denoted simply by $S\times T$. 
Given a scheme $S$, we denote by $p_S:  S \to \Spec(k)$ the structure map. 

\sssec{}

By default, we only consider schemes of finite type, and we simply call them \virg{schemes}. We denote by $\Aff$ the category of affine schemes of finite type over $k$. We do not use derived algebraic geometry in this paper. 
Let $\bbP^{\infty}$  be the infinite projective space $\bbP^\infty := \colim_{d \gg 0} \bbP^d$, viewed as an ind-scheme.

\sssec{}

By an \emph{$\infty$-category} we mean an $\left(\infty,1\right)$-category, for example, in the sense of Joyal and Lurie \cite{Lurie}.
We denote by $\Spc$ the $\infty$-category of $\infty$-groupoids (that is, topological spaces). In this paper, by the word \emph{prestack} we mean a functor $\Y: \Aff^\op \to \Spc$.

\sssec{}

We will work with DG categories (regarded as stable $\infty$-categories tensored over $k$).
For instance, we denote by $\Vect_k$ the DG category of
chain complexes of vector spaces over $k$, modulo quasi-isomorphism (whose homotopy
category is equivalent to the derived category of the ordinary category of $k$-vector spaces).
Since the field $k$ is mostly fixed, we often adopt the notation $\Vect = \Vect_k$.
We abuse notation and refer to objects of $\Vect$ simply as \virg{vector spaces}.
\sssec{}

Similarly, $\QCoh(\Y)$ and $\Dmod(\Y)$ will denote the DG categories of quasi-coherent sheaves and D-modules on a prestack $\Y$. See e.g. \cite{DG} for a foundational treatment. We have natural identifications
\begin{gather}
\nonumber
\Dmod(\Spec(k)) \simeq  \Vect_k
\\
\nonumber
\QCoh(\Spec(k)) \simeq  \Vect_k.
\end{gather}

\sssec{}

The functorialities of $\QCoh$ and $\Dmod$ are spelled out for instance in \cite{DG}.
At some point, we will need to use the DG category $\ICoh$ of ind-coherent sheaves on a formal completion. We refer to \cite{Book} for the functorialities of $\ICoh$ and for the related notion of \emph{nil-isomorphism}.

\ssec{Contractibility} \label{ssec:prelim-contract}

\sssec{}

Let $\Y$ be a scheme, a pseudo-scheme (see below) or an algebraic stack. We define the de Rham cohomology $H^*(\Y)$ of a prestack $\Y$ via D-modules as follows:
$$
H^*(\Y) := \CHom_{\Dmod(\Y)}(\omega_\Y, \omega_\Y).
$$
Its homology is defined as 
$$
H_*(\Y) := (p_\Y)_!(\omega_\Y),
$$
where $(p_\Y)_!$ is the partially-defined left adjoint to $(p_\Y)^!$. It turns out that the value of $(p_\Y)_!$ on $\omega_\Y$ is well-defined: this is because $\omega_\Y$ is holonomic.
By adjunction, we see that $H_*(\Y)$ is pre-dual to $H^*(\Y)$.

\medskip

We say that $\Y$ is \emph{homologically contractible} if $H_*(\Y) \simeq k$; equivalently, if $H^*(\Y) \simeq k$.
Occasionally, we abuse terminology and just say that $\Y$ is \emph{contractible}.

\sssec{}

The above definitions are justified by the Riemann-Hilbert correspondence, see for instance the discussion in \cite[Section 6.1]{Barlev}. 
Namely, in case there is an embedding $\iota:  k \hto \bbC$,
we have an isomorphism
$$
H_*(\Y) \otimes_k \bbC 
\simeq
 H^{\Top}_*(\Y^{\Top}, \bbC) 
 $$
of chain complexes of $\bbC$-vector spaces, where $H^{\Top}_*$ is the ordinary singular homology and $\Y^{\Top}$ the topological space (in the classical analytic topology) associated to $\Y$.

\sssec{}

A \emph{pseudo-scheme} is a prestack that can be written as a colimit (not necessarily filtered) of schemes along proper maps. When $\Y = \colim_a Y_a$ is a pseudo-scheme, its DG category of D-modules is well-behaved: it can be written as a colimit 
$$
\Dmod(\Y) \simeq \colim_{a \in \CA} \Dmod(Y_a)
$$
along the $!$-pushforwards $(f_{a \to b})_!: \Dmod(Y_a) \to \Dmod(Y_b)$.
In particular, $\Dmod(\Y)$ is compactly generated and we have 
$$
H_*(\Y)
\simeq
\colim_{a \in \CA} H_*(Y_a).
$$

\begin{rem}
Suppose that $j: \CX \to \Y$ is a locally closed embedding of prestacks and $\Y$ is a pseudo-scheme. Then $\CX$ is a pseudo-scheme, too.
\end{rem}

\sssec{}

When $\Y = \colim_a Y_a$ is a pseudo-scheme, its \emph{Borel-Moore homology} is well-defined as
$$
\Hren(\Y) := (p_\Y)_{*,\dR}(\omega_\Y),
$$
see \cite{omega} for a detailed discussion. Concretely, using the presentation, we obtain
$$
\Hren(\Y)
\simeq
\colim_{a \in \CA} \Hren(Y_a).
$$
We say that $\Y$ is pseudo-proper if there is a presentation $\Y = \colim_{a} Y_a$ with all the $Y_a$ proper. In this case, $\Hren(\Y) \simeq H_*(\Y)$.

\sssec{} \label{sssec:pseudo-proper map}

More generally, we say that a map $f: \CX \to \Y$ of prestacks is \emph{pseudo-proper} if, for any scheme $S \to \Y$, the fiber-product $\CX \times_\Y S$ can be presented as a pseudo-scheme $\colim_a Z_a$, with each $Z_a$ proper over $S$.
If $f:\CX \to \Y$ is pseudo-proper, the functor $f_!: \Dmod(\CX) \to \Dmod(\Y)$ is well-defined and satisfies base-change against $!$-pullbacks. It is easy to check that pseudo-properness is preserved by composition and base-change.

\begin{lem}[{\cite[Lemma 6.1.7]{AG2}}] \label{lem:k'-fibers}
For $f:\CX \to \Y$ pseudo-proper, the natural arrow $f_!(\omega_{\CX}) \to \omega_\Y$ is an isomorphism iff, for any field $k' \supseteq k$ and any $y \in \Y(k')$, the $k'$-pseudo-scheme $\CX_y := \CX \times_\Y y$ is homologically contractible (that is $H_*(\CX_y) \simeq k'$).
\end{lem}

\sssec{} \label{sssec:ICoh-contr}

Let $\Y$ be a prestack. We say that $\Y$ is $\O$-contractible if $H^*(\Y, \O) \simeq k$. We shall consider this only for $\Y = \colim_{a} Y_a$ a pseudo-scheme, in which case
$$
H^*(\Y,\O) \simeq \lim_{a \in A} H^*(Y_a, \O).
$$
In case $\Y$ is a pseudo-proper pseudo-scheme, we define also the notion of $\ICoh$-contractibility. This means that the natural map
$$
H_*(\Y,\O) := p_{\Y,*}^{\ICoh}(\omega_{\Y}^{\ICoh}) \longto k
$$
is an isomorphism. Since $\Y$ is pseudo-proper, by adjunction, 
$
H_*(\Y,\O) \simeq \colim_a H_*(Y_a, \O)
$
is pre-dual to 
$$
\CHom_{\ICoh(\Y)}
\bigt{
\omega_{\Y}^{\ICoh},
\omega_{\Y}^{\ICoh}
}.
$$
The latter is isomorphic to
$
\CHom_{\QCoh(\Y)}
\bigt{
\O_\Y,
\O_\Y
}
\simeq
H^*(\Y,\O)
$
by means of the fully-faithful functor $\Upsilon_{\Y}: \QCoh(\Y) \to \ICoh(\Y)$.

\sssec{}

An important example of a pseudo-proper pseudo-scheme is the \emph{Ran space} of a smooth projective curve $X$, see \cite{BD-chiral, contract}. Let $\fSet$ be the 1-category of nonempty finite sets and surjective maps. Consider the functor 
$$
\fSet^\op \longto \mathrm{Schemes}
$$
$$
I \squigto X^I
$$
$$
\phi=[I \tto J] \squigto \Delta_{\phi}: X^J \hto X^I.
$$
Then
$$
\Ran := 
\colim_{I \in \fSet^\op} X^I.
$$
This pseudo-scheme is homologically contractible (as proven in \cite{BD-chiral} and recalled in \cite[Section 6]{contract}) and $\O$-contractible (the same proof works).
The operation of union of finite sets of $X$ endows $\Ran$ with the structure of (non-unital) commutative monoid.

\sssec{} \label{sssec:Ran-indscheme}

When $\Y$ is a prestack equipped with a map $\Y \to \Ran$, we usually denote by $\Y_{X^I}$ or by $\Y_I$ its restriction $\Y \times_{\Ran} {X^I}$. We say that $\Y \to \Ran$ is a \emph{$\Ran$-indscheme} if each $\Y_{X^I}$ is an indscheme.

\ssec{Domains, generic maps and generic sections}

\sssec{}

Given an affine $k$-scheme $S$, regarded as a test scheme, we let $X_S:=S \times X$. Following \cite{Barlev}, we say that an open subset $U \subset X_S$ is a \emph{domain}
if the induced map $U \hto X_S \tto S$ is surjective (equivalently, if all the fibers of $U \to S$
are dense in $X$).
To indicate that $U$ is a domain, we use the notation
$U \Subset X_S$.

More generally, let $U \Subset X_S$ be a domain and $U' \subset U $ an open subset. We say that $U' \subset U$ is a \emph{sub-domain} (and use the notation $U' \Subset U$) if $U' \subseteq X_S$ is also a domain. 

\sssec{}

Following \cite{Barlev}, we denote by $\Dom_X$ the following $1$-category. Its objects are pairs $(S,U)$, with $S$ an affine scheme and $U \subseteq X_S$ a domain. An arrow $(S,U) \to (T,V)$ in $\Dom_X$ is a map $S \to T$ in $\Aff$ such that the $U \hto X_S \to X_T$ factors through $V$.

In this paper, we will often deal with functors out of $\Dom_X$ and their left Kan extensions along the forgetful functor $\Dom_X \to \Aff$ that sends $(S,U) \squigto S$.

\sssec{}

For instance, letting $\Y$ be a prestack, consider the functor
$$
\Maps(X, \Y)^{\Dom_X}:
\Dom_X \to \Spc, 
\hspace{.4cm}
(S,U) 
\squigto 
\CHom_{\PreStk}(U,\Y),
$$
where $\CHom_{\PreStk}$ denotes the Hom-spaces in the $\infty$-category $\PreStk$.
The left Kan extension of $\Maps(X, \Y)^{\Dom_X}$ along $\Dom_X \to \Aff$  yields the prestack of generic maps $X \dasharrow \Y$:
$$
\GMaps(X,\Y):
S \squigto
\colim_{U \Subset X_S} 
\CHom_{\PreStk}(U,\Y). 
$$
Now suppose that $\Y = Y$ is a scheme. Then, by \cite[Section 3]{Barlev}, $\GMaps(X,Y)$ is a pseudo-scheme and it is pseudo-proper if the $Y$ is proper.

\sssec{}

For $\pi: \CZ \to X$ a prestack over $X$, we define the prestack $\GSect(X,\CZ)$ of generic sections of $\pi$ in the same manner. 
As a functor out of $\Dom_X$, it sends $(S, U \subseteq X_S)$ to the space of maps $U \to \CZ$ compatible with $X$.
Then we set
$$
\GSect(X,\CZ)(S) 
=
\colim_{U \Subset X_S} \Sect(X,\CZ)^{\Dom_X}(S,U).
$$
When $\CZ$ is a (proper) $X$-scheme, $\GSect(X,\CZ)$ is a (pseudo-proper) pseudo-scheme.

\sssec{} \label{sssec:quasi-sections}

Related to the notions of generic maps/sections are the notion of quasi-maps and quasi-sections, see \cite{Barlev}. We will only need the definition of $\QSect(X, \bbP E)$ for a vector bundle $E \to X$. An $S$-point of $\QSect(X, \bbP E)$ is given by $(L, i: L \hto E)$, where $L$ is a line bundle on $X_S$ and $i$ an injection of coherent sheaves with $S$-flat cokernel.
It turns out that $\QSect(X, \bbP E)$ is a scheme with a natural map to $\GSect(X, \bbP E)$.

\sssec{}

Let $Y$ be a scheme  and let $F \subset E$ be an inclusion of (algebraic) vector bundles on $Y$. We usually do not distinguish between vector bundles and their locally free sheaves of sections. 
We say that $F$ is a \emph{subbundle} of $E$ if the injection $F \hto E$ of coherent sheaves has locally free cokernel. We use the notation $F \prec E$ to describe this situation.

\sssec{}

Let $Z$ be a scheme and $E$ be a vector bundle on $Z$. When there is no risk of confusion, for a map $f:  Y \to Z$ of schemes, we denote by $\restr E Y$ the pullback vector bundle $f^*(E)$. A similar notation will be used when restricting $G$-bundles, or $G$-local systems.

We stress the fact that the notation $\restr E Y$ is used in cases where $f$ is not an immersion. For example, we will use it for the maps $U \hto X_S \to X$ and $U \hto X_S \tto S$, with $U$ a domain.

\ssec{Opers}

We review the notion of $G$-oper following \cite{BD-opers}.

\sssec{}
  
Recall that $G$ always denotes a connected reductive group over $k$. We fixed a Borel subgroup $B \subset  G$ once and for all, with unipotent radical $N$. 
Denote by $\Par$ the poset of standard (relative to $B$) parabolic subgroups of $G$. Let $\I$ be the Dynkin diagram of $G$ and, for $P \in \Par$, set $\I_P \subseteq \I$ to be the Dynkin diagram of the Levi quotient of $P$.
For brevity, define $\Par' := \Par \smallsetminus \{G\}$.

\sssec{}

We also fix a Cartan subgroup $T \subset  B$.
Let $\g \supset {\fb} \supset {\ft}$ be the corresponding Lie algebras.
Denote by $\sS \subset {\mathfrak \ft}^*$ the set of simple roots with respect to ${\mathfrak b}$.
There is a unique Lie algebra grading (the \emph{principal grading})
$\g = \bigoplus_m \g^{(m)}$ such that $\g^{(0)} = {\mathfrak \ft}$,
$\g^{(1)} =\bigoplus_{\alpha \in \sS} \g^\alpha$ and
$\g^{(-1)} = \bigoplus_{\alpha \in \sS} \g^{-\alpha}$. 
 The corresponding filtration 
\begin{equation} \label{eqn:filtration-Lie-alge}
\g^{(\geq r)} = \bigoplus_{m \geq r} \g^{(m)},
\mbox{    for $r \in \bbZ$}
\end{equation}
is $B$-invariant.

\sssec{}

\nc{\cF}{\CF}
\nc{\calT}{\mathcal{A}}

Let $S$ be an affine scheme and $p: U \to S$ a smooth map of dimension $1$. In practice, $U$ is a domain in $X_S$ and $p: U \hto X_S \to S$ is the projection map.
Let $H$ be an algebraic group with Lie algebra $\fh$ and let  $\pi: E \to U$ an $H$-bundle on $U$. Since $p$ and $\pi$ are both smooth, we have a short exact sequence
\begin{equation} \label{eqn:rel-tang-sequence}
0 \to T_{E/U} \to T_{E/S} \to \pi^*(T_{U/S}) \to 0
\end{equation}
of $H$-equivariant vector bundles on $E$.
Thus, this sequence corresponds to an exact sequence of vector bundles on $U = E/H$: this is the Atiyah exact sequence 
\begin{equation} \label{eqn:atiyah}
0
\to
\fh_E
\to
\calT_E
\to
T_{U/S}
\to 0.
\end{equation}
Here
$$
\calT_E := (T_{E/S})/H \to E/H = U
$$
is a vector bundle on $U$ of rank $\dim H + 1$, while $\fh_E$ denotes the $E$-twist of the adjoint $H$-representation $\fh$.

\sssec{}

Now let $E$ be a $G$-bundle on $U$, equipped with a $B$-reduction $F$. Thus, $F$ is a $B$-bundle on $U$ together with an identification $E := G \times^B F$.
The Atiyah exact sequences for $E$ and $F$ combine into a commutative diagram
\begin{equation} 
\nonumber
\begin{tikzpicture}[scale=1.5]
\node (00) at (0,0) {$ \fg_F = \fg_E$};
\node (10) at (1.5,0) {$ \calT_E$};
\node (01) at (0,1) {$ \fb_F$};
\node (11) at (1.5,1) {$\calT_F$};
\node (20) at (3,0) {$ T_{U/S}$.};
\node (21) at (3,1) {$T_{U/S}$};
\path[right hook ->,font=\scriptsize,>=angle 90]
(00.east) edge node[above] {$ $}  (10.west); 
\path[right hook->,font=\scriptsize,>=angle 90]
(01.east) edge node[above] {$ $} (11.west); 
\path[->>,font=\scriptsize,>=angle 90]
(10.east) edge node[above] {$ $}  (20.west); 
\path[->>,font=\scriptsize,>=angle 90]
(11.east) edge node[above] {$ $} (21.west); 
\path[right hook ->,font=\scriptsize,>=angle 90]
(01.south) edge node[right] {$ $} (00.north);
\path[right hook->,font=\scriptsize,>=angle 90]
(11.south) edge node[right] {$ $} (10.north);
\path[->,font=\scriptsize,>=angle 90]
(21.south) edge node[right] {$=$} (20.north);
\end{tikzpicture}
\end{equation}

\sssec{}

Since the filtration \eqref{eqn:filtration-Lie-alge} is $B$-invariant, it induces a filtration $\g^{(\geq r)}_{F}$ of $\g_F$ by sub-bundles.
We obtain a non-negative filtration of ${\calT}_{E}$ by sub-bundles as follows: for $r \geq 0$, set $\calT^{(\geq r)}_{E} \subset \calT_{E}$ to be the preimage of 
$$
\g^{(\geq r)}_F / {\fb}_F
\subset 
\g_F / {\fb}_F
\xto{\simeq}
 \calT_{E}/\calT_F
$$
along the projection $  \calT_{E} \tto  \calT_{E}/\calT_F$.
In particular, we have $\calT_E^{\geq 0} =\calT_F$ while
$$
\calT_E^{(\geq -1)}
/
\calT_E^{(\geq 0)} 
=
\g^{(\geq -1)}_F / \g^{(\geq 0)}_F
\simeq 
\bigoplus_{\alpha \in \sS} \g^{-\alpha}_F.
$$
Here $\g^{-\alpha}_F$ is the $F$-twist of the $B$-representation $\g^{-\alpha}$, the $B$-action on $\g^{-\alpha}$ being defined via the composition $B \tto T$.

\sssec{}

Let $U/S$ be as above and $E$ a $G$-bundle on $U$.
By definition, a \emph{connection on $E$ relative to $S$} is a right splitting
$$
\nabla: T_{U/S} \to \calT_E
$$
of Atiyah's exact sequence \eqref{eqn:atiyah}. 
If a $E$ is also equipped with a $B$-reduction $F$, the connection and the $B$-reduction combine to givean $\O_U$-linear map
$$
T_{U/S} \xto{\nabla} \calT_E\tto \calT_E/ \calT_{F}
\simeq \g_F /\fb_F,
$$
that is, a section $\nabla_{/F} \in H^0(U, \g_F / \fb_F \otimes \Omega_{U/S})$.

\sssec{}
Given a $G$-local system $\sigma = (E,\nabla)$ on $U/S$, an \emph{extended oper structure} on $\sigma$ is a $B$-reduction $F$ of $E$ subject to the condition:
\begin{itemize}
\item 
$\nabla(T_{U/S}) \subset \calT_E^{(\geq -1)}
\subset \calT_{E}$.
\end{itemize}
An \emph{oper structure} on $\sigma$ is an extended oper structure $F$ subject to the following further condition:
\begin{itemize}
\item 
for each simple root $\alpha \in \sS$, the composition
$$
T_{U/S} \xrightarrow{\;\; \nabla\;\;} 
{\calT_E^{(\geq -1)}/\calT_F} 
\tto \g^{-\alpha}_F  
$$
is an isomorphism of line bundles.
\end{itemize}

\ssec{Generic opers}\label{ssec:defn-gen-opers}

We now define $\Opbargen_G$, the prestack classifying $G$-local systems equipped with a generic extended oper structure.

\sssec{}

Recall the $1$-category $\Dom_X$ of pairs $(S,U)$, with $S \in \Aff$ and $U \subseteq X_S$ a domain. We first define a functor
$$
\Opbar_G^{\Dom_X}:
\Dom_X \longto \Spc.
$$ 
This functor assigns to $(S,U)$ the groupoid of triples $(E,\nabla, F)$, where:
\begin{itemize}
\item
$E$ is a $G$-bundle on $X_S$;
\item
$\nabla: T_{X_S/S} \to \calT_E$ is an $S$-relative connection on $E$;
\item
$F$ is an extended oper structure on $(\restr E U, \restr \nabla U)$.
\end{itemize}

\sssec{}

We then define the prestack $\Opbargen_{G}$ to be the assignment
$$
\Opbargen_{G} (S)
:=   
\colim_{U \Subset S \times X} 
\Opbar_G^{\Dom_X}{(S,U)}.
$$
Hence, $\Opbargen_{G,\sigma}$ assigns to a test scheme $S \in \Aff$ the groupoid of equivalence classes of the following data:
\begin{itemize}
\item
$(E, \nabla)$ as above (that is, an $S$-family of $G$-local systems on $X$);
\item
a domain $U \Subset X_S$;
\item
an extended oper structure on $\restr \sigma U$.
\end{itemize} 

\begin{rem}
Similarly, we define $\Opgen_{G}$, the prestack of generic opers, to be as above but with the word \virg{extended} omitted. The inclusion $\Opgen_{G} \hto \Opbargen_{G}$ is an open embedding.
\end{rem}

\begin{lem} \label{lem:pibar-proper}
The forgetful map $\ol\pi: \Opbargen_G \to \LSG$ is pseudo-proper.
\end{lem}

\begin{proof}
By definition, $\ol\pi$ factors as 
$$
\Opbargen_G
\hto
\LS_G \times_{\Bun_G} \Bun_G^{B\ggen}
\to
\LSG,
$$
with the left map a closed embedding. Hence, we only need to show that the second map is pseudo-proper. By base-change, it suffices to prove this for the forgetful map $\Bun_G^{B\ggen} \to \Bun_G$. The latter is proven in \cite[Section 4]{Barlev} using Drinfeld's compactification.
\end{proof}

\sssec{}

We now look at the fibers of the above forgetful map.
Let $\sigma$ be a $k$-point of $\LSG$, that is, a $G$-bundle with connection $(E,\nabla)$ on $X$. We set 
$$
\Opbargen_{G,\sigma} := \Opbar_G \utimes_{\LSG} \sigma,
\hspace{.4cm}
\Opgen_{G,\sigma} := \Op_G \utimes_{\LSG} \sigma.
$$
The above lemma implies that these are pseudo-schemes, with the former pseudo-proper.

\sec{Reduction to the case of extended opers} \label{sec:Steinberg}

In this section, we prove Theorem \ref{mainthm:St-Opers} assuming Theorem \ref{mainthm:contract}. We then show that the latter is equivalent to Theorem \ref{mainthm:contract-fiber}.

\ssec{Constructing a canonical map}

Recall that $\pi: \Opgen_G \to \LSG$ denotes the structure forgetful map. We first construct a natural arrow
$$
\xi_G: \pi_{*,\dR} (\omega_{\Opgen_G})
\longto
\St_G
$$
in $\Dmod(\LSG)$. Then we will use Theorem \ref{mainthm:contract} to prove that $\xi_G$ is an isomorphism.

\sssec{}

By Lemma \ref{lem:pibar-proper}, the structure map
$$
\ol\pi: \Opbargen_G 
\longto \LSG
$$
is pseudo-proper, hence $\ol\pi_{*,\dR} \simeq \ol\pi_!$. Thus, adjunction yields an arrow
$$
\ol\pi_{*,\dR} 
\bigt{
\omega_{\Opbargen_G}
}
\xto{\ol\xi_G}
\omega_{\LSG}
$$
which is an isomorphism by Theorem \ref{mainthm:contract}.

\sssec{}

The obvious map $j: \Opgen_G \hto \Opbargen_G$ is an open embedding.\footnote{In passing, recall that $j$ is an isomorphism over the locus of irreducible $G$-local systems.}
The complementary closed subspace can be described as follows.
For each standard parabolic $P \in \Par$, consider the prestack $\Opbargen_{G,P}$ classifying degenerate $G$-opers that are zero on the nodes of $\I \setminus \I_P$. Denote by $i_P: \Opbargen_{G,P} \hto \Opbargen_G$ the obvious closed embedding. 
The prestack
$$
\CZ
:=
\colim_{P \in \Par'}
\;
\Opbargen_{G,P}
$$
is closed inside $\Opbargen_G$ and it is complementary to $\Opgen_G$. 
It follows that 
$$
j_{*,\dR} (\omega_{\Opgen_G})
\simeq
\coker
\Bigt{
\colim_{P \in \Par'} \;
(i_P)_{*,\dR}
\bigt{
\omega_{\Opbargen_{G,P}}
}
\to
\omega_{\Opbargen_G}
}
$$
and thus that 
\begin{eqnarray} \label{eqn:writing St via opers}
\pi_{*,\dR} (\omega_{\Opgen_G})
& \simeq &
\ol\pi_{*,\dR} \circ
j_{*,\dR} (\omega_{\Opgen_G})
\\
\nonumber
& \simeq &
\coker
\Bigt{
\colim_{P \in \Par'} \;
\ol\pi_{*,\dR} \circ
(i_P)_{*,\dR}
\bigt{
\omega_{\Opbargen_{G,P}}
}
\longto
\ol\pi_{*,\dR} \bigt{
\omega_{\Opbargen_G}}
}.
\end{eqnarray}

\sssec{}

Observe now that, for each $P$, there is a pseudo-proper forgetful map
$$
\ol\pi_P: 
\Opbargen_{G,P}
\to
\LSP
$$
and consequently a natural arrow
$$
(\ol\pi_P)_{*,\dR} 
\bigt{
\omega_{\Opbargen_{G,P}}
}
\xto{\ol\xi_P}
\omega_{\LSP}.
$$

\sssec{}

For each inclusion $Q \subseteq P$ in $\Par$, the square 
$$ 
\begin{tikzpicture}[scale=1.5]
\node (01) at (0,1) {$\Opbargen_{G,Q}$};
\node (11) at (2,1) {$\Opbargen_{G,P}$};
\node (00) at (0,0) {$\LS_Q$};
\node (10) at (2,0) {$\LSP$};
\path[->,font=\scriptsize,>=angle 90]
(00.east) edge node[above]  {$ $}  node[below] { } (10.west);
\path[->,font=\scriptsize,>=angle 90]
(01.south) edge node[left] {$\ol\pi_Q$}  (00.north);
\path[->,font=\scriptsize,>=angle 90]
(01.east) edge node[above] {$ $ } node[below] { } (11.west);
\path[->,font=\scriptsize,>=angle 90]
(11.south) edge node[right] {$\ol\pi_P$} node[left] {} (10.north);
\end{tikzpicture}
$$
is evidently commutative. More precisely, the maps $\ol\pi_P$ assemble into a natural transformation $\ol\pi_?: \Opbargen_{G,?} \to \LS_?$.
In particular, by applying this to $P \subseteq G$, we obtain isomorphisms
$$
\ol\pi_{*,\dR}
(i_P)_{*,\dR}
\bigt{
\omega_{\Opbargen_{G,P}}
}
\simeq
(\fp_P)_{*,\dR}
\Bigt{
(\ol\pi_P)_{*,\dR}
\bigt{
\omega_{\Opbargen_{G,P}}
}},
$$
which are obviously compatible with the induction maps as $P$ varies.\footnote{Recall that $\fp_P: \LS_P \to \LS_G$ denotes the induction map.}
Then the maps $\ol\xi_P$ can be applied together to the RHS of \eqref{eqn:writing St via opers} to yield an arrow 
$$
\xi_G: \pi_{*,\dR}(\Opgen_G) \to \St_G,
$$ 
which is the morphism we were looking for.

\ssec{Contractibility statements}

The above discussion implies that $\xi_G$ is an isomorphism if the same is true for the maps $\ol\xi_P$'s for all $P \in\Par$: in other words, Theorem \ref{mainthm:St-Opers} would follow if we can prove that, for each $P \in \Par$, adjunction gives the isomorphism
$$
(\ol\pi_P)_{*,\dR}
\bigt{
\omega_{\Opbargen_{G,P}}
}
\simeq
\omega_{\LSP}.
$$
For $P=G$, this is exactly the statement of Theorem \ref{mainthm:contract}. In the general case, the assertion follows from Theorem \ref{mainthm:contract} applied to the  Levi quotient $M$ of $P$, which is still a classical reductive group. Indeed, it suffices to observe that the diagram
$$ 
\begin{tikzpicture}[scale=1.5]
\node (01) at (0,1) {$\Opbargen_{G,P}$};
\node (11) at (2,1) {$\Opbargen_M$};
\node (00) at (0,0) {$\LS_P$};
\node (10) at (2,0) {$\LS_M$};
\path[->,font=\scriptsize,>=angle 90]
(00.east) edge node[above]  {$ $}  node[below] { } (10.west);
\path[->,font=\scriptsize,>=angle 90]
(01.south) edge node[left] {$\ol\pi_P$}  (00.north);
\path[->,font=\scriptsize,>=angle 90]
(01.east) edge node[above] {$ $ } node[below] { } (11.west);
\path[->,font=\scriptsize,>=angle 90]
(11.south) edge node[right] {$\ol\pi_M$} node[left] {} (10.north);
\end{tikzpicture}
$$
is cartesian (this is immediate from the definition of $\Opbargen_{G,P}$) and then apply base-change.

\sssec{}

Let us now show that Theorem \ref{mainthm:contract-fiber} and Theorem \ref{mainthm:contract} are equivalent. The implication from the former to the latter is base-change. For the opposite implication, since $\ol\pi$ is pseudo-proper, we can use Lemma \ref{lem:k'-fibers} and reduce to checking the homological contractibility of the $k'$-fibers.
However, it is easy to see that $\LSG \otimes _k k'$ is the $k'$-stack of $G_{k'}$-local systems on $X_{k'}$ (and similarly for $\Opbar_G$). Hence, up to renaming $k'$ with $k$, we are back to the statement of Theorem \ref{mainthm:contract-fiber}.

\sssec{} \label{ssec:example-GL2}

The remainder of this paper (with the exception of Section \ref{sec:Tsen}) will be devoted to prove Theorem \ref{mainthm:contract-fiber}: for $\sigma = (E, \nabla) \in \LSG(k)$ arbitrary, $H_*(\Opbargen_{G,\sigma}) \simeq k$. 
Here we give the simplest example.
Let $G = GL_2$ and $\sigma = (E, \nabla)$ a $G$-local system on $X$ (that is, a vector bundle of rank $2$ with a flat connection).
In this case, the closed embedding $\Opbargen_{G, \sigma} \hto \fL_E$, with $\fL_E$ the space of generic line sub-bundles of $E$, is obviously an isomorphism.
Hence, the theorem is a consequence of the contractibility of $\fL_E$, which is well-known (\cite{Barlev}) since $\fL_E \simeq \GMaps(X, \bbP^1)$.

\sec{Extended opers for classical groups} \label{sec:classical}

From now on, we consider only the case where the group $G$ is classical and begin our proof of Theorem \ref{mainthm:contract-fiber}. The goal of this section is to reduce the contractibility of $\Opbargen_{G,\sigma}$ to that of a simpler space $\Y_{G,\sigma}$. The definition of the latter does depend on the Dynkin type of $G$, so we need to treat each type separately (except for type $B$ and $D$, which can go together).

\ssec{General strategy}

Let us describe the general outline of our construction.

\sssec{}

Since $G$ is classical, $G$-local systems and extended $G$-opers can be easily reinterpreted as more explicit structures. 
As a preliminary observation, pairs of groups $G$ and $G'$ of the same semisimple type have isomororphic flag varieties and thus isomorphic spaces of extended oper structures. 
Hence, it suffices to treat general linear groups, symplectic groups and orthogonal groups.

\sssec{}

In each of these cases, we have a standard $G$-representation $V$ and a maximal parabolic $P$ such that $G/P \subseteq \bbP V$.
Our first step consists of writing down a closed embedding 
$$
\Opbargen_{G,\sigma}
\hto
\GMaps(X, G/B).
$$
We will then consider the string
$$
\phi:
\Opbargen_{G,\sigma}
\hto
\GMaps(X, G/B)
\to
\GMaps(X, G/P)
\hto
\GMaps(X, \bbP V),
$$
the middle arrow being induced by the projection $G/B \tto G/P$. This composition is far from being surjective, in general. We will be able to identify the image of $\phi$ with a certain explicit closed subspace $\Y_{G,\sigma} \subseteq \GMaps(X, \bbP V)$. 

\sssec{}

We then show that the resulting map
$
\xi: 
\Opbargen_{G,\sigma}
\to
\Y_{G,\sigma}
$
has contractible fibers. We check this by induction: we will see that these fibers are isomorphic to spaces of generic extended opers for reductive groups of the same type as $G$, but of strictly smaller semisimple rank.
To carry out this computation, we rely on the construction of Section \ref{ssec:amplify bundles}, which explains how the combination of a line bundle and a linear differential operator generates partial oper structures.

\sssec{}

The contractibility of the fibers of $\xi$ implies that $\Opbargen_{G,\sigma}
$ and $\Y_{G,\sigma}$ have the same de Rham cohomology. Now, $\Y_{G,\sigma}$ is a closed subspace of $\GMaps(X, \bbP V)$ defined by \virg{not too many} explicit equations. Counting the number of these equations and their degrees will let us prove the contractibility of $\Y_{G,\sigma}$ in Section \ref{sec:Oka}, after the somewhat lengthy preparation of Section \ref{sec:contract-LE}.

\ssec{A key construction} \label{ssec:amplify bundles}

\sssec{}

Let $X'$ a smooth (not necessarily proper) curve and $E$ a vector bundle on $X'$. 
Let $D: E \to E$ be a linear differential operator and $F \prec E$ a sub-bundle. In practice, $D$ will be the operator $\nabla_\nu$ corresponding to a connection $\nabla$ on $E$ and $\nu \in T_{X'}$.

\sssec{}

Observe that the induced map $\wt D: F \to E/F$ is $\O_{X'}$-linear. This follows from the definition of \emph{linear} differential operator. For instance, in the case of $D = \nabla_\nu$, for $\phi \in \O_{X'}$ and $v \in F$, we have
$$
\nabla_{\nu}( \phi \cdot v) 
= 
{\nu}(\phi ) \cdot v + \phi \cdot \nabla_{\nu}(v) 
\equiv
 \phi \cdot \nabla_{\nu}(v)  \,\, (\mathrm{mod} \, F).
$$

\sssec{}

We define $\langle F,D(F) \rangle \in \Coh(X')^\heart$ to be the coherent sheaf 
$$
\langle F, D(F) \rangle
:=
\ker 
\Bigt{
E \to (E/F)/\wt D(F) 
}.
$$
Since $E$ a vector bundle and $X$ is a smooth curve, $\langle F, D(F) \rangle$ is a vector bundle too.

\sssec{}

We have tautological injections
$$
F \hto \langle F, D(F) \rangle \hto E
$$
of coherent sheaves. We claim that $F \hto \langle F, D(F) \rangle$ is a bundle embedding, that is, its cokernel is locally free. Indeed, observe that 
$$
\coker
\Bigt{
F 
\to 
\ker 
\bigt{
E \to (E/F)/\wt D(F) 
}
}
\simeq
\wt D(F).
$$
The latter is a coherent sub-sheaf of the locally free sheaf $E/F$, hence it is locally free since $X'$ has dimension one.

\sssec{}

Up to replacing $X'$ with a nonempty open sub-variety, we can assume that the inclusion $\langle F, D(F) \rangle \hto E$ is a bundle embedding. We can thus repeat the above construction: the sub-bundle $F_1: = \langle F, D(F) \rangle \prec E$ yields a new vector bundle 
$$
\langle F, D(F), D^2(F) \rangle
:=
\langle F_1, D(F_1) \rangle
$$
sitting in a flag
$$
F \prec
\langle F, D(F) \rangle
\prec
\langle F, D(F), D^2(F) \rangle 
\hto 
E.
$$
Obviously, this process stabilizes after a finite number of iterations, yielding a sub-bundle that is preserved by the differential operator $D$.

\ssec{General linear groups}

In this section, we prove Theorem \ref{mainthm:contract-fiber} in the case of $G = GL_n$. Let $\sigma = (E, \nabla)$ be a $G$-local system. We will construct a map 
$$
\xi: \Opbargen_{G,\sigma} 
\to
 \GMaps(X,\bbP^{n-1})
$$
that is pseudo-proper and with the following property: the arrow 
$$
\xi_!(\omega_{\Opbargen_{G,\sigma} }) \to \omega_{ \GMaps(X,\bbP^{n-1})}
$$ 
is an isomorphism of D-modules. 
Using this, the contractibility of $\Opbargen_{G,\sigma}$ boils down to the well-known contractibility of $\GMaps(X,\bbP^{n-1})$. A new proof of the latter fact, along with generalizations, will be provided in Section \ref{sec:contract-LE}.

\sssec{}

A $G$-local system on $X$ is just a pair $(E, \nabla)$ where $E$ is an $n$-dimensional vector bundle on $X$ and 
$
\nabla: E \longto E \otimes \Omega_X 
$
is a map satisfying the Leibniz rule.
Thus the functor $\Opbar_{G,\sigma}^{\Dom_X}$ sends $(S,U) \in \Dom_X$ to the set 
$\Opbar^{\Dom_X}_{G,\sigma}(S,U)$ consisting of those full flags of sub-bundles
$$
0 = F_0 \prec F_1 \prec \cdots \prec F_{n-1} \prec  F_n = \restr E U
$$
that satisfy the following condition:
\begin{itemize}
\item
for $0\leq i \leq n-2$, the map $\restr \nabla U$ restricts to a map
$
F_i \to F_{i+1} \otimes \restr{\Omega_X}U.
$
\end{itemize}
As usual, we have:
$$
\Opbargen_{G,\sigma} (S) 
=
\colim_{U \Subset X_S} \Opbar^{\Dom_X}_{G,\sigma}(S,U).
$$

\sssec{} \label{sssec:simplifying-assumptions}

It turns our that we can simplify our above description of $\Opbargen_{G,\sigma}$ even further.
Let $X' \subseteq X$ be a nonempty open affine sub-variety where $E$ and $T_X$ are both trivial, with trivializations $E \simeq \O_{X'}^{\oplus n}$ and $T_{X'} \simeq \O_{X'}$ (the latter given by a vector field that we call $\nu$). Let $k[X']:= H^0(X', \O_{X'})$. In this case, the connection is determined by the linear differential operator
$$
\nabla_\nu = \nu + A: 
 \O_{X'}^{\oplus n}
  \to
   \O_{X'}^{\oplus n}
$$
for some matrix $A \in \mathrm{Mat}_{n,n}(k[X'])$.

\sssec{}

The functor
$$
\Dom_{X} \to \Dom_{X'}
$$
$$
(S,U) \squigto (S, U \cap X')
$$
is evidently cofinal. Hence, it suffices to define our $\Opbar_{G,\sigma}$ on $\Dom_{X'}$. In this case, we have the following description.
The value of $\Opbar_{G,\sigma}^{\Dom_{X'}}$ on $(S,U) \in \Dom_{X'}$ is the set of full flags of subbundles
\begin{equation} \label{eqn:flag-for-GLn-simple}
0 = F_0 \prec F_1 \prec \cdots \prec F_{n-1} \prec  \O_U^{\oplus n}
\end{equation}
that satisfy the following condition:
\begin{itemize}
\item
for $0\leq i \leq n-2$, the map $\restr{\nabla_{\nu}} U: \O_U^{\oplus n} \to \O_U^{\oplus n}$ restricts to a map
$F_i \to F_{i+1}$.
\end{itemize}

\sssec{}

Consider now the natural transformation
$$
\xi^{\Dom_{X'}}:
\Opbar_{G,\sigma}^{\Dom_{X'}}
 \to
 \Maps^{\Dom_{X'}}(X, \bbP^{n-1} )
$$
whose value on $(S,U) \in \Dom_{X'}$ is the arrow that sends the sequence \eqref{eqn:flag-for-GLn-simple} to $F_1 \prec \O_U^{\oplus n}$.
Passing to colimits, we obtain a map
$$
\xi: 
\Opbargen_{G,\sigma}
 \to
\GMaps(X, \bbP^{n-1}).
$$

\begin{lem}
The map $\xi$ is pseudo-proper.
\end{lem}

\begin{proof}
Let $P \subseteq G$ be the maximal parabolic corresponding to $(n-1,1)$-blocks. Obviously, $\bbP^{n-1} \simeq G/P$, and thus $\xi$ factors as 
$$
\Opbargen_{G,\sigma}
\hto
\GMaps(X, G/B)
\to
\GMaps(X, G/P).
$$
Since the left map is a closed embedding, we just need to show that the right map is pseudo-proper. 
This is a particular case of the fact that $\Bun_G^{B \ggen} \to \Bun_G^{P \ggen}$ is pseudo-proper.
\end{proof}

\begin{prop}
The natural arrow 
$$
\xi_!(\omega_{\Opbargen_{G,\sigma}})
\to 
\omega_{\GMaps(X, \bbP^{n-1})}
$$
is an isomorphism in $\Dmod( \GMaps(X, \bbP^{n-1}))$.
\end{prop}

\begin{proof}
We proceed by induction on $n$, the cases $n=1,2$ being obvious (see Section \ref{ssec:example-GL2} for the latter).
Since $\xi$ is pseudo-proper, by Lemma \ref{lem:k'-fibers} we can check the assertion at the level of $k'$-fibers, for all fields $k' \supseteq k$. Namely, for a $k'$-point $\ell: \Spec(k') 
\to 
\GMaps(X, \bbP^{n-1} )$, we need to show that the pseudo-scheme
$$
 \Opbargen_{G,\sigma}
 \utimes_
 {
\GMaps(X, \bbP^{n-1} ) 
 }
 \Spec(k') 
$$
is homologically contractible over $k'$.
It is easy to realize that, up to renaming, we can assume $k' =k$.
Thus, $\ell$ can be represented by $(U, L_U \prec \O_U^{\oplus n})$ for some nonempty $U \subseteq X'$. We now show that to $\ell$ we can canonically attach an integer $1 \leq d \leq n$ and a $GL_d$-local system defined on some nonempty open $U' \subset U$.
Indeed, by applying the construction of Section \ref{ssec:amplify bundles} repeatedly, we can find a certain nonempty open $U' \subset U$ and a sub-bundle of $\restr E{U'} = \O_{U'}^{\oplus n }$ given by
$$
F:=
\langle L_{U'}, \nabla_\nu(L_{U'}), \ldots, \nabla^{n-1}_\nu(L_{U'}) \rangle.
$$
Letting $d$ denote the rank of $F$, we see that
$$
F \simeq 
\langle L_{U'}, \nabla_\nu(L_{U'}), \ldots, \nabla^{d-1}_\nu(L_{U'}) \rangle
$$
and that $F$ is preserved by $\nabla_\nu$. Thus, the pair $\sigma_\ell:=(F, \nabla_\nu)$ is a $GL_d$-local system defined on $U'$. 
Let $\wt\sigma_\ell = (E/F, \wt\nabla)$ denote the resulting $GL(n-d)$-local system. 

Since $\sigma_\ell$ comes with an extended oper structure given by the flag $(F_i)_{0 \leq i  \leq d}$ with 
$$
F_i = 
\langle L_{U'}, \nabla_\nu(L_{U'}), \ldots, \nabla^{i-1}_\nu(L_{U'}) \rangle,
$$
it is clear that 
$$
 \Opbargen_{G,\sigma}
 \utimes_
 {
 \GMaps(X, \bbP^{n-1} )
 }
 \ell 
\simeq
\Opbargen_{GL(n-d), \wt\sigma_\ell}.
$$
We then apply the induction hypothesis and conclude the proof.
\end{proof}

\sssec{}

From here, the proof of Theorem \ref{mainthm:contract-fiber} is obvious: the above proposition ensures that
$$
H_*(\Opbargen_{G,\sigma})
\to
H_*(\GMaps(X,\bbP^{n-1}))
$$
is an isomorphism, and the latter equals $k$. As we see below, the situation is more interesting for the other classical groups.

\sssec{} \label{sssec:generic-flags}

In preparations for these cases, let us fix the following notation. For $1 \leq n \leq p$, consider the parabolic subgroup $Q \subseteq GL_p$ that parametrizes flags of type 
$$
W_1 \subset \cdots \subset W_n \subset k^{\oplus p}
$$
with $\dim(W_i)=i$ for $1 \leq i \leq n$.
We denote by 
$$
\GFlags_{p,n} := \GMaps(X, GL_p/Q)
$$ 
the resulting space of generic flags of sub-bundles of $\O_X^{\oplus p}$.

\ssec{Symplectic groups} \label{ssec:type-Cn}

Now let us consider the case where $G = Sp_{2n}$.

\sssec{}

An $Sp_{2n}$-local system is a triple $\sigma =(E, \nabla, \omega)$, where $(E,\nabla)$ is a $GL_{2n}$-local system on $X$ and $\omega$ is  
a non-degenerate alternating form
$$
\omega: \Lambda^{2} E \to \CO_X
$$ 
that is horizontal with respect to $\nabla$: for $\mu \in T_X$, and $e_1, e_2 \in E$, we have
$$
\mu(\omega(e_1,e_2)) = 
\omega( 
\nabla_{\mu}(e_1),e_2
 )  
 + 
\omega( e_1,\nabla_{\mu}(e_2) ).
$$ 

\begin{rem} \label{rem:Sp-trick with nabla}
Let $\mu \in T_X$ and $e \in E$. The above formula implies that
\begin{equation} \label{eqn:nabla-trick-Sp(2n)}
\mu(\omega(e, \nabla_{\mu}(e))) = 
\omega( 
\nabla_{\mu}(e),\nabla_{\mu}(e)
 )  
 + 
\omega( e,\nabla_{\mu}^2 (e) )
=
\omega( e,\nabla_{\mu}^2 (e) ).
\end{equation}
This formula will be useful later on.
\end{rem}

\sssec{}

The flag variety of $Sp_{2n}$ parametrizes full flags of $\omega$-isotropic spaces. Hence, unraveling the definitions, for $(S,U) \in \Dom_X$, the set $\Opbar_{G,\sigma}^{\Dom_X}(S,U)$
consists of those sequences of sub-bundles
$$
0 = F_0 \prec F_1 \prec \cdots \prec F_{n} \prec \restr E U
$$
that satisfy the following conditions:
\begin{itemize}
\item 
for $0\leq i \leq n$, we have
$\mathrm{rank}(F_i)= i$;
\item
$F_n$ is Langrangian with respect to the form $\omega$, i.e. $\omega(F_n, F_n) = 0$;
\item 
for $0\leq i \leq n-1$, the map $\restr \nabla U$ restricts to a map
\begin{equation} 
\label{eqn:nabla-U for Sp(2n)}
\restr \nabla U :  F_i \to F_{i+1} \otimes \restr{\Omega_X}U.
\end{equation} 
\end{itemize}

\sssec{}

As above, we work generically and thus can restrict to a nonempty affine open $X' \subset X$ such that $\restr E{X'} \simeq \O_{X'}^{\oplus 2n}$ trivial, and $T_{X'}$ is trivial (generated by a derivation $\nu$).
In particular, $\nabla$ is determined by the differential operator $\nabla_\nu = \nu + A: E \to E$, with $A \in \mathfrak{sp_{2n}}(k[X'])$. We see that $\Opbargen_{G,\sigma}$ is the functor
$$
 S 
\mapsto 
\left\{
[U, F_1 \prec \cdots \prec F_{n} \prec \O_U^{\oplus 2n}] \in 
\GFlags^{\leq n}_{2n}(S)
\;
\bigg| 
\;
\begin{array}{l}
\omega(F_{n}, F_{n}) = 0  \\
\nabla_\nu(F_i) \subseteq F_{i+1} 
\text{ for $1 \leq i \leq n-1$} 
 \end{array}
\right\}.
$$
Our auxiliary space $\Y_{G,\sigma}$ is defined by the assignment
$$
S
\mapsto
\left\{
[U, L \prec \O_U^{\oplus 2n}] \in 
\GMaps(X,\bbP^{2n-1})(S)
\;
\bigg| 
\;
\omega(\nabla^i_\nu(L), \nabla_\nu^j(L)) = 0 
\text{ for $0 \leq i < j \leq n-1$}
\right\}.
$$

\begin{warning}
For a single pair $(i,j) \neq (0,0)$, the equation
$$
\omega(\nabla_\nu^{i}(L), \nabla_\nu^{j}(L)) = 0
$$
is ill-defined. However, a simple application of the Leibniz rule shows that the totality of the equations (for all $0 \leq i,j \leq n-1$) yields well-defined conditions on $L$.
\end{warning}

\sssec{}

There is an evident pseudo-proper map 
$$
\xi: \Opbargen_{G, \sigma} \longto \Y_{G,\sigma}
$$ 
that sends a generic flag to its first member. To check that $\xi$ does indeed land in $\Y_{G,\sigma}$, we just need to verify that, for $[U, F_\bullet] \in \Opbargen_{G, \sigma}(S)$ and $L =F_1$, the generic vector bundle
$$
F:= 
\langle
L, 
\nabla_\nu(L), 
\ldots,
\nabla^{n-1}_\nu(L)
\rangle
$$
is isotropic: this is certainly true since $F$ is a generic sub-bundle of $F_n$, and the latter is Langrangian.

\begin{prop}
The map $\xi$ has contractible fibers.
\end{prop}

\begin{proof}
Using Lemma \ref{lem:k'-fibers} and the usual argument, we just need to check the contractibility of the $k$-fibers. Thus, let $[U,L] \in \Y_{G,\sigma}(k)$ and let $1 \leq d \leq n$ be the rank of the generic vector bundle $F$ defined above.
We deduce
$$
F 
\simeq 
\langle
L, 
\nabla_\nu(L), 
\ldots,
\nabla^{d-1}_\nu(L)
\rangle
$$
is an isotropic sub-bundle $F \prec E$ preserved by $\nabla_\nu$. It follows that $\nabla$ descends to a connection on $F^\perp/F$, thereby definining an $H$-local system $\wt \sigma$, with $H= Sp_{2(n-d)}$, on a nonempty open subset of $U$.

By construction, we see that the fiber of $\xi$ above $[U,L]$ is isomorphic to $\Opbargen_{H, \wt\sigma}$. The latter is homologically contractible by induction.
\end{proof}

\sssec{}

It turns out that many of the equations defining $\Y_{G,\sigma}$ inside $\GMaps(X,\bbP^{2n-1})$ are redundant. Namely, we just need the following ones:
$$
\omega(\nabla^i_\nu(L), \nabla_\nu^{i+1}(L)) = 0 
\text{ for $0 \leq i \leq n-2$}.
$$
In fact, we obtain the remaining ones by applying Remark \ref{rem:trick-SO} several times.
All in all, Theorem \ref{mainthm:contract-fiber} for $G=Sp_{2n}$ amounts to proving the homological contractibility of the pseudo-scheme $\Y_{G,\sigma}$ given by
$$
S
\mapsto
\left\{
[U, L \prec \O_U^{\oplus 2n}] \in 
\GMaps(X,\bbP^{2n-1})(S)
\;
\bigg| 
\;
\omega(\nabla^i_\nu(L), \nabla_\nu^{i+1}(L)) = 0 
\text{ for $0 \leq i \leq n-2$}
\right\}.
$$
Thus, $\Y_{G,\sigma}$ is defined by $n-1$ quadratic equations in $\GMaps(X,\bbP^{2n-1})$. We will see that the inequality $2n-1 \geq 2(n-1)$ is the key condition for the contractibility of $\Y_{G,\sigma}$.

\ssec{Orthogonal groups}

Finally, let us consider the orthogonal groups $SO_{2n-1}$ and $SO_{2n}$. These can be treated uniformly: we let $G=SO_m$ with $m$ denoting either $2n-1$ or $2n$, so that $n = \lceil m/2 \rceil$.

\sssec{}

An $SO_m$-local system is a triple $\sigma =(E, \nabla, \beta)$, where $(E,\nabla)$ is a $GL_m$-local system on $X$ and $\beta$ is  
a non-degenerate symmetric bilinear form
$$
\beta: \Sym^{2} E \to \CO_X
$$ 
with the following property: for any $\mu \in T_X$, and $e_1, e_2 \in E$, we have
$$
\mu(\beta(e_1,e_2)) = 
\beta( 
\nabla_{\mu}(e_1),e_2
 )  
 + 
\beta( e_1,\nabla_{\mu}(e_2) ).
$$ 
There is also a horizontal trivialization of $\det(E)$, which we omit because it is irrelevant in what follows.

\begin{rem} \label{rem:trick-SO}
A particular case of the above formula is
$$
\mu(\beta(e,e)) = 
\beta( 
\nabla_{\mu}(e),e
 )  
 + 
\beta( e,\nabla_{\mu}(e) )
=
2 \beta( e,\nabla_{\mu}(e) ).
$$ 
This will be useful later.
\end{rem}

\sssec{}

The flag variety of $G = SO_m$ parametrizes partial flags of $\beta$-isotropic spaces of dimensions $1 \leq i \leq n-1$. We emphasize that this holds true for both $SO_{2n}$ and $SO_{2n-1}$. 

Unraveling the definition of Griffiths transversality in the present case, we see that the set $\Opbar^{\Dom_X}_{G,\sigma}(S,U)$ consists of those sequences of sub-bundles
$$
0 = F_0 \prec F_1 \prec \cdots \prec F_{n-1} \prec \restr E U
$$
that satisfy the following conditions:
\begin{itemize}
\item 
for $0\leq i \leq n-1$, we have
$\mathrm{rank}(F_i)= i$;
\item
$F_{n-1}$ is $\beta$-isotropic, i.e. $\beta(F_{n-1}, F_{n-1}) = 0$;
\item 
for $0\leq i \leq n-2$, the map $\restr \nabla U$ restricts to a map
$$ 
\restr \nabla U :  F_i \to F_{i+1} \otimes \restr {\Omega_X} U.
$$
\end{itemize} 

\sssec{}

As above, we work generically and thus assume that $X$ is affine, $E \simeq \O_X^{\oplus m}$ is trivial, and $T_X$ is trivial (generated by a derivation $\nu$).
In passing, we remairk that $\nabla_\nu = \nu + A$ for some $A \in \mathfrak{so}_m(k[X'])$.

We then see that $\Opbargen_{G,\sigma}$ is the functor
$$
S 
\mapsto 
\left\{
[U, F_1 \prec \cdots \prec F_{n-1} \prec \O_U^{\oplus m}] \in 
\GFlags^{\leq n-1}_{m}(S)
\;
\bigg| 
\;
\begin{array}{l}
\beta(F_{n-1}, F_{n-1}) = 0  \\
\nabla_\nu(F_i) \subseteq F_{i+1} 
\text{ for $0 \leq i \leq n-2$} 
 \end{array}
\right\}.
$$
Our auxiliary space $\Y_{G,\sigma}$ is defined by the assignment
$$
S
\mapsto
\left\{
[U, L \prec \O_U^{\oplus m}] \in 
\GMaps(X,\bbP^{m-1})(S)
\;
\bigg| 
\;
\beta(\nabla^i_\nu(L), \nabla_\nu^j(L)) = 0 
\text{ for $0 \leq i \leq j \leq n-2$}
\right\}.
$$

\sssec{}

There is an evident pseudo-proper map $\xi: \Opbargen_{G, \sigma} \to \Y_{G,\sigma}$ which sends a generic flag to its first member. Indeed, the conditions defining $\Y_{G,\sigma}$ mean that the generic vector bundle
$$
F:= 
\langle
L, 
\nabla_\nu(L), 
\ldots,
\nabla^{n-2}_\nu(L)
\rangle
$$
is $\beta$-isotropic. For $[U, F_\bullet] \in \Opbargen_{G, \sigma}$ and $L =F_1$, this is certainly true as $F$ is a generic sub-bundle of $F_{n-1}$.

\begin{prop}
The map $\xi$ has contractible fibers.
\end{prop}

\begin{proof}
As before, we just need to check the contractibility of the $k$-fibers. Thus, let $[U,L] \in \Y_{G,\sigma}(k)$ and let $1 \leq d \leq n-1$ be the rank of the generic vector bundle $F$.
We deduce
$$
F 
\simeq 
\langle
L, 
\nabla_\nu(L), 
\ldots,
\nabla^{d-1}_\nu(L)
\rangle
$$
and that $F \prec E$ is preserved by $\nabla_\nu$. It follows that $\nabla$ descends to a connection on $F^\perp/F$, thereby definining an $H$-local system $\wt \sigma$, with $H= SO_{m-2d}$, on a nonempty open subset of $U$.

By construction, we see that the fiber of $\xi$ above $[U,L]$ is isomorphic to $\Opbargen_{H, \wt\sigma}$. The latter is homologically contractible by induction.
\end{proof}

\sssec{}

It turns out that many of the equations defining $\Y_{G,\sigma}$ inside $\GMaps(X,\bbP^{m-1})$ are redundant. Namely, we just need the following ones:
$$
\beta(\nabla^i_\nu(L), \nabla_\nu^i(L)) = 0 
\text{ for $0 \leq i \leq n-2$}.
$$
In fact, we obtain the remaining ones by applying Remark \ref{rem:trick-SO} several times.

All in all, Theorem \ref{mainthm:contract-fiber} for $G=SO_m$ amounts to proving the homological contractibility of the pseudo-scheme $\Y_{G,\sigma}$ defined by the assignment
$$
S
\mapsto
\left\{
[U, L \prec \O_U^{\oplus m}] \in 
\GMaps(X,\bbP^{m-1})(S)
\;
\bigg| 
\;
\beta(\nabla^i_\nu(L), \nabla_\nu^i(L)) = 0 
\text{ for $0 \leq i \leq \lceil m/2 \rceil -2$}
\right\}.
$$
Thus, $\Y_{G,\sigma}$ is defined by $(\lceil m/2 \rceil -1)$ quadratic equations in $\GMaps(X,\bbP^{m-1})$. We will see that the inequality $m-1 \geq 2(\lceil m/2 \rceil -1) $ is the key condition for the contractibility of $\Y_{G,\sigma}$.

\sec{Contractibility of the space of generic lines} \label{sec:contract-LE}

In the previous section, we reduced the proof of Theorem \ref{mainthm:contract-fiber} to the homological contractibility of an auxiliary space $\Y_{G,\sigma}$.
In the case where $G$ is of type $A_n$, this auxiliary space is simply $\fL:= \GMaps(X, \bbP^n)$, the space of generic line sub-bundles of $\O_X^{\oplus (n+1)}$.
Slightly more generally, we will consider the space 
$$
\fL_E := \GSect(X, \bbP E)
$$ 
of generic line sub-bundles of a fixed vector bundle $E \to X$. Our goal in this section is to prove the $\O$-contractibility and the homological contractibility of $\fL_E$. 
To achieve this goal, we will write $\fL_E$ as a quotient (in the fppf topology) of two explicit $\Ran$-indschemes\footnote{Recall this definition and the related notation from Section \ref{sssec:Ran-indscheme}.} whose cohomologies we can control.

\ssec{Ran-families of infinite projective spaces}

Let $E \to X$ be a vector bundle on $X$ and denote by
$\fL_E := 
\GSect(X, \bbP E)$
the space of generic sections of $\bbP E \to X$. Thus, $\fL_E$ parametrizes generically defined line sub-bundles of $E$.

\begin{rem} \label{rem:can-trivialize-E}
As already observed and used earlier, removing any finite number of points from $X$ does not affect the spaces $\GMaps(X, Y)$ or $\GSect(X,Y)$. In particular, restricting to a nonempty open subset $X' \subseteq X$ where $E$ is trivial, we see that 
$$
\GSect(X, \bbP E)
\simeq 
\GMaps(X', \bbP^{r-1})
\simeq
\GMaps(X, \bbP^{r-1})
$$
with $r =\mathrm{rank}(E)$. However, having $E$ not necessarily trivial makes our constructions below neater. 
\end{rem}

\sssec{}

We wish to present $\fL_E$ using the Ran space. To this end, we will define a $\Ran$-indscheme $\fL_{E,\Ran}$ equipped with a natural map to $\fL_E$.
The idea behind $\fL_{E,\Ran}$ is very natural: we look for a space that parametrizes pairs $(\ul x, \ell)$, where $\ul x \in \Ran$ and $\ell$ a point of the infinite-dimensional projective space $\bbP \bigt{ H^0(X \smallsetminus D_{\ul x}, E)}$.
A formal definition requires setting up some notation.

\sssec{}

Denote by 
$$
D^{\inc} \subseteq \Ran \times X \to \Ran
$$
the incidence relative divisor, by $U^{\inc}$ its complement and by $\pi^{\inc}_{\Ran}: U^{\inc} \to \Ran$ the projection.
Similarly, for a finite set $I \in \fSet$, let $D^{\inc}_I \subseteq X^I \times X$ be the incidence divisor, $U^{\inc}_I$ its complement and $\pi^{\inc }_{X^I}$ the projection onto $X^I$. We also denote by $\pi_{\Ran}: \Ran \times X \tto \Ran$ and by $ \pi_{X^I}: X^I \times X \tto X^I$ the projections, that is, the relative compactifications of $\pi^{\inc }_{\Ran}$ and $\pi^{\inc }_{X^I}$.

The following lemma is a standard consequence of Grauert theorem and Riemann-Roch.

\begin{lem} \label{lem:locally free}
For any $I \in \fSet$ there is $d_0 \geq 0$ such that 
$$
( \pi_{X^I})_*
\Bigt{
\!\!
\restr E {X^I \times X}(d D^{\inc}_I)
}
\in \Coh(X^I)^{\heart}
$$ 
is locally free whenever $d \geq d_0$.
\end{lem}

\sssec{} 

We now set
$$
\fL_{E,\Ran}
:=
\bbP_{\Ran} 
\Bigt{
(\pi^{\inc }_{\Ran})_*( \restr  E {U^{\inc}})
}.
$$
Tautologically, we have
$$
\fL_{E,\Ran}
\simeq
\colim_{I \in \fSet^\op}
\fL_{E,I},
$$
with
$$
\fL_{E,I}
:=
\bbP_{X^I} 
\Bigt{
(\pi_{X^I}^{\inc})_*( \restr  E {U_I^{\inc}})
}
\simeq
\colim_{d \gg 0}
\;
\bbP_{X^I} 
\Bigt{
( \pi_{X^I})_*
\bigt{
\restr E {X^I \times X}(d D^{\inc}_I)
}
}.
$$
For later use, we also set
$$
\fL_{E,I}^{\leq d}
:=
\bbP_{X^I} 
\Bigt{
(\pi_{X^I})_*
\bigt{
\restr E {X^I \times X}(d D^{\inc}_I)
}
}.
$$
Thus, each $\fL_{E,I}$ is an ind-proper indscheme and $\fL_{E,\Ran}$ is a $\Ran$-indscheme ((and, in particular, a pseudo-proper pseudo-scheme).

Informally, we view $\fL_{E,\Ran}$ as a family of $\bbP^\infty$'s over $\Ran$. Since $\bbP^\infty$ and $\Ran$ are both $\O$-contractible, we expect the same to be true for $\fL_{E,\Ran}$. Indeed:

\begin{prop} \label{prop:fL-E-Ran contractible}
The pseudo-proper $\Ran$-indscheme $\fL_{E,\Ran}$ is $\ICoh$-contractible and $\O$-contractible.
\end{prop}

\begin{proof}
As explained in Section \ref{sssec:ICoh-contr}, $\O$-contractibilty and $\ICoh$-contractibility are equivalent. Let us prove the latter.
Fix $I \in \fSet$, denote by 
$$
p_I: \fL_{E,I} \to X^I
$$
$$
p_I^{\leq d}: \fL_{E,I}^{\leq d} \to X^I
$$
the structure maps. Let us first show that, for fised $I$ and $d$ large enough as in the lemma above, the natural map
$$
(p_I^{\leq d})_*^{\ICoh}
\Bigt{
\omega_{\fL_{E,I}^{\leq d}}
}
 \longto 
 \omega_{X^I}
$$
is an isomorphism. This can be checked Zariski-locally on $X^I$: taking an open subset $U \subseteq X^I$ on which $(\pi_{X^I})_*
\bigt{
E_{X^I \times X}(d D^{\inc}_I)
} \in \Coh(X^I)$ is free, we are in the product situation of $\bbP^N \times U \to U$ and our claim is just the $\ICoh$-contractibility of $\bbP^N$.

Now, the isomorphism
$$
(p_I)_*^{\ICoh}(\omega_{\fL_{E,I}})
\simeq
\colim_{d \gg 0 } \;
(p_I^{\leq d})_*^{\ICoh}
\Bigt{
\omega_{\fL_{E,I}^{\leq d}}
}
$$
implies that the natural map
$$
(p_I)_*^{\ICoh}
\bigt{
\omega_{\fL_{E,I}}
}
 \longto 
 \omega_{X^I}
$$
is an isomorphism, too. It follows that
$$
\Gamma^{\ICoh}(\fL_{E,\Ran}, \omega)
\simeq
\Gamma^{\ICoh}(\Ran, \omega).
$$
To conclude, we observe that the well-known proof of the homological contractibility of $\Ran$, see \cite{BD-chiral, contract}, can be repeated verbatim to yield $\Gamma^{\ICoh}(\Ran, \omega) \simeq k$.
\end{proof}

\ssec{A presentation of $\fL_E$} \label{ssec:FE-construction}

Let $E \to X$ be a vector bundle as above. In this section, we construct a natural map
$F_E:
\fL_{E, \Ran} \to \fL_E$.
This maps will turn out to be a surjection in the fppf topology.

\sssec{}
Constructing a map $F_E: \fL_{E, \Ran} \to \fL_E$ amounts to providing a system of maps
$$
F_{E,I}^{d}:
\fL_{E, I}^{\leq d} \longto \fL_E.
$$
that are compatible with $I$ and $d$ in the natural way. %
We will actually construct a compatible system of maps
$$
\QF_{E,I}^d:
\fL_{E,X^I}^{\leq d} \to \QSect(X, \bbP E)
$$
and thus a map
$$
\QF_{E} : \fL_{E,\Ran} \to  \QSect(X, \bbP E).
$$
Here $\QSect(X, \bbP E)$ denotes the scheme of \emph{quasi-sections} of $\bbP E \to X$, see Section \ref{sssec:quasi-sections}.
Our map $F_E$ will then be the composition
$$
\fL_{E,\Ran}
\xto{\QF_E}
\QSect(X, \bbP E)
\longto
\GSect(X, \bbP E) = \fL_E.
$$

\sssec{} \label{sssec:S-point of fL-I-d}

An $S$-point of $\fL_{E, I}^{\leq d}$ consists of a pair $q=(\ul x, \eta)$ where:
\begin{itemize}
\item
$\ul x = (x_i)$ is a $I$-tuple of elements of $X(S)$; we denote by $D_{\ul x} \subseteq X_S$ the associated divisor and by $U_{\ul x} \Subset X_S$ its complement;
\item
$\eta$ is the datum of a line sub-bundle on $S$:
$$
\eta: L_S 
\prec
 (\pi_S)_*( \restr E {X_S} (d D_{\ul x})),
$$
where $\pi_S: X_S \to S$ is the projection.

\end{itemize}

\sssec{}

By adjunction,  $\eta$ gives rise to an arrow
$$
 (\ol\pi_S)^*(L_S) 
\to 
\restr E {X_S} (d D_{\ul x}),
$$
and thus to an arrow
$$
\eta':
 (\ol\pi_S)^*(L_S) (-d D_{\ul x})
\to 
\restr E {X_S}.
$$

\begin{lem} \label{lem:S-flatness}
The above map $\eta'$ is an injection of coherent sheaves with $S$-flat cokernel.
In other words, $\eta'$ is an $S$-point of $\QSect(X, \bbP E)$.
\end{lem}

\begin{proof}
The question is Zariski-local on $S$; to simplify the situation, let us assume that $L_S = \O_S$ so that $\eta \in H^0(X_S, \restr E {X_S} (d D_{\ul x}))$.
Now the flatness of $\coker(\eta)$ amounts to 
$$
\eta_s \in H^0(X_s, E (d D_{\ul x,s}))
$$
being nonzero for all closed points $s \in S$.
This condition immediately implies that $\eta'$ is injective as a map of coherent sheaves on $X_S$.
Then the $S$-flatness of the cokernel of $\eta'$ is equivalent to $\eta'_s: \O_X(- d D_{\ul x, S}) \to E$ being injective for all closed points $s \in S$. Since $X$ is a smooth curve, this follows again from the non-vanishing of $\eta_s$.
\end{proof}

\sssec{}

The assignment $q= (\ul x, \eta) \squigto \eta'$ defines our map $\QF_{E,I}^d(S) \to \QSect(X, \bbP E)(S)$. It is straightforward to verify that such a construction is functorial in $I$, $d$ and $S$. We thus obtain the map $\QF_E$.
To obtain $F_E$, we post-compose with $\QSect(X, \bbP E) \to \GMaps(X, \bbP E)$. Concretely, the latter map sends $\eta'$ as above to the equivalence class of $[U, \restr {\eta'} U]$, where $U \subseteq X_S$ is the maximal open locus where the $\coker(\eta')$ is $U$-flat. We know by \cite[Lemma 3.2.8]{Barlev} that $U$ is a domain of $X_S$.

\begin{example}\label{example:coord}
Let $q =(\ul x, \eta)$ be an $S$-point of $\fL_{E,\Ran}$ as in Section \ref{sssec:S-point of fL-I-d}.
Assume furthermore that $L_S \simeq \O_S$ (this is of course always true Zariski-locally on $S$) and that $E \to X$ is the trivial vector bundle of rank $n$. In this case, $F_E(q)$ can be described in coordinates.
Namely, for some $d \in \bbN$ sufficiently large, $\eta$ is the datum of an $n$-tuple $(f_i)$ of sections in $H^0(X_S, \O_{X_S}(d D_{\ul x}))$, with the following property: for any closed point $s \in S$, at least one of the $f_i(s) \in H^0(X-D_{x_i,s}, \O)$ is non-zero.
A different isomorphism $L_S \simeq \O_S$ changes each $f_i(s)$ by a common scalar from $k^\times$. In particular, $\Gamma_q \subset X_S$, the zero locus of the $f_i$'s, is well-defined and, by the above property, the complementary open subset
$$
U_q := (X_S \smallsetminus D_{\ul x}) \smallsetminus \Gamma_q
$$ 
is a domain. Now consider the inclusion
$$
i_q:\O_{U_q} \hookrightarrow \O_{U_q}^{\oplus n}
 $$
defined by the restriction of the section $\eta$ to $U_q$. It is clear that $i_q$ is the inclusion of a line sub-bundle with $U_q$-flat cokernel. By construction, $F_{E}(q) = [U_q, i_q]$. Note that changing all the $f_i$'s by an invertible element of $H^0(S,\O_S)$ does not change the image of $i_q$.
\end{example}

\sssec{}

Recall from \cite[Proposition 3.2.2]{Barlev} that the natural map $\QSect(X, \bbP E) \to \GSect(X, \bbP E)$ is an effective epimorphism in the Zariski topology. Here we show that a similar fact holds true for $F_{E}$. This will eventually allow us to write $\fL_E$ as the quotient of $\fL_{E,\Ran}$ by a monoid object in $\Ran$-indschemes.

\begin{prop} \label{prop:FE-epi}
The map $F_{E} : \fL_{E,\Ran} \to \fL_E$ constructed above is an effective epimorphism for the fppf topology.
\end{prop}

\begin{proof}

Let $S$ be an affine scheme and $\ell \in \fL_E(S)$.
By definition, $\ell$ is represented by $[U,L_U]$, where $U \Subset X_S = S \times X$ is a domain and $L_U$ a $1$-dimensional sub-bundle $L_U \prec \restr E U$.
 We wish to show that
there exist
\begin{itemize}
\item
an fppf cover $\wt S \to S$
\item
a finite set $I$ and an integer $d \geq 0$
\item
an $\wt S$-point $q$ of $\fL_{E,I}^{\leq d}$
\end{itemize}
such that $F_{E,I}^{\leq d}(q) = \restr \ell {\wt S}$.
We will proceed in steps.

\sssec*{Step 1}

Here, we need to use some specifics from the proof of \cite[Lemma 3.2.9]{Barlev}, which states that the map $\QSect(X, \bbP E) \to \GSect(X, \bbP E)$ is an effective epimorphism in the Zariski topology.
The lemma states that there exists an arrow $(S',U') \to (S,U)$ in $\Dom_X$ such that:
\begin{itemize}
\item
the resulting map $S' \to S$ is a Zariski cover;
\item
the pullback of $L_U$ to $U'$ is trivial;
\item
$U'$ is a divisor complement, i.e. there exists a line bundle $Q$ on $ X_{S'}$, equipped with a nonzero section $\sigma$, such that $U' = X_{S'} \smallsetminus Z(\sigma)$;

\item

for some $d \geq 0$, there is an injective map
\begin{equation} \label{eqn:kappa map}
\O_{ X_{S'}} \hto \restr E {X_{S'}} \otimes Q^{\otimes d}
\end{equation} 
of coherent sheaves, with $S'$-flat cokernel, whose restriction to some sub-domain of $U'$ is isomorphic to $\restr \ell {S'}$.
\end{itemize}

\sssec*{Step 2}

Let 
\begin{equation} \label{eqn:kappa prime}
\nonumber
\kappa': \O_{ X_{S'}} \otimes Q^{\otimes (-d)}
\hto \restr E {X_{S'}} 
\end{equation}
be the map induced from \eqref{eqn:kappa map} by twisting.
This way, $\kappa'$ is an $S'$-point of $\QSect(X, \bbP E)$.
It remains to find a fppf cover $S'' \to S'$ and an element $q \in \fL_{E, \Ran}(S'')$ such that $\QF_E(q)$ equals $\restr {\kappa'}{S''}$.
Choose an appropriate fppf cover $S'' \to S'$ such that the support of the divisor $\restr \sigma {S''}$ can be written as the union of the graphs of some maps $x_i: S'' \to X$.
In other words, we obtain an isomorphism $\restr Q {S''} \simeq \O_{X_{S''}}(D_{\ul x})$ for some point $\ul x \in  X^I(S'')$.
Hence, $\kappa'' := \restr {\kappa'}{S''}$ reads as
$$
\kappa'': \O_{ X_{S''}} (-d D_{\ul x})
\hto
\restr E {X_{S''}} .
$$

\sssec*{Step 3}

By adjunction and twisting, $\kappa''$ yields a map
$$
\kappa:  \O_{S''}
\to
(\pi_{S''})_*
\bigt{
\restr E {X_{S''}} (d D_{\ul x})
}.
$$ 
We claim that $\kappa$ is an injection of coherent sheaves on $S''$ with flat cokernel. Indeed, as in the proof of Lemma \ref{lem:S-flatness}, both properties follow easily from the corresponding properties of $\kappa''$.

Hence $q := (\ul x, \kappa)$ determines an $S''$-point of $\fL_{E,X^I}^{\leq d}$. It is immediate from the definition of $\QF_{E,I}^d $ that $\QF_{E,I}^d (q) = \kappa''$, as desired.
\end{proof} 

\ssec{Compatibility with products}

Consider the obvious map
$$
\mu:
\fL_{E_1} \times \fL_{E_2} \to \fL_{E_1 \otimes E_2}.
$$
In this section, we will construct an analogous map
$$
\mu_{\Ran}:
\fL_{E_1,\Ran}
\times 
\fL_{E_2,\Ran} 
\to 
\fL_{E_1 \otimes E_2,\Ran}
$$
and exploit it to define a natural monoid object in $\Ran$-indschemes.

\sssec{}

Let $(\ul x, \eta)$ and $(\ul y, \theta)$ be $S$-points of $\fL_{E_1,X^I}$ and $\fL_{E_2, X^J}$. Explicitly, 
$\eta$ and $\theta$ are given by line sub-bundles
$$
\eta: L_S 
\prec
 (\pi_S)_*
 \bigt{ \!\!
 \restr {E_1}{X_S} (d D_{\ul x})
 }
$$
$$
\theta: M_S 
\prec
 (\pi_S)_*
 \bigt{\!\!
  \restr {E_2}{X_S} (e D_{\ul y})
  }
$$
for some $d,e \gg 0$.
By tensoring these together and using the lax monoidal structure on $\ol\pi_*$, we obtain a map
$$
L_S \otimes M_S 
\to
 (\pi_S)_*
 \bigt{ \!\!
 \restr {E_1}{X_S} (d D_{\ul x})
 }
 \otimes
 (\pi_S)_*
 \bigt{\!\!
  \restr {E_2}{X_S} (e D_{\ul y})
  }
 \to
  (\pi_S)_*
  \Bigt{\!\!
  \restr {E_1 \otimes E_2}{X_S} 
  (d D_{\ul x} + eD_{ \ul y})
  }.
$$
It is easy to see that this arrow in $\Coh(S)$ is injective with flat cokernel, hence the same is true for its post-composition with the injection
$$
  (\pi_S)_*
  \Bigt{\!\!
  \restr {E_1 \otimes E_2}{X_S} 
  (d D_{\ul x} + eD_{ \ul y})
  }
  \hto
    (\pi_S)_*
  \Bigt{\!\!
  \restr {E_1 \otimes E_2}{X_S} 
  ((d+e) D_{\ul x \cup \ul y})
  }.
$$
We thus obtain an $S$-point of $\fL_{E_1 \otimes E_2, \Ran}$, more precisely, an $S$-point of  $\fL_{E_1 \otimes E_2, X^{I \sqcup J}}^{\leq (d+e)}$. This procedure is functorial in $S$, thereby yielding a map $\mu_{\Ran}$ as desired.

\sssec{}

The map $F_E: \fL_{E,\Ran} \to \fL_E$ constructed in Section \ref{ssec:FE-construction} is compatible with products in the following sense.

\begin{lem} \label{lem:comm-diagram-of-projectivizations}
The diagram 
\begin{equation} 
\nonumber
\begin{tikzpicture}[scale=1.5]
\node (00) at (0,0) {$ \fL_{E_1 \otimes E_2, \Ran}$};
\node (10) at (4,0) {$ \fL_{E_1 \otimes E_2}$};
\node (01) at (0,1) {$ \fL_{E_1, \Ran} \otimes \fL_{E_2, \Ran}$};
\node (11) at (4,1) {$\fL_{E_1} \times \fL_{E_2}$};
\path[->,font=\scriptsize,>=angle 90]
(00.east) edge node[above] {$F_{E_1 \otimes E_2}$}  (10.west); 
\path[->,font=\scriptsize,>=angle 90]
(01.east) edge node[above] {$F_{E_1} \times F_{E_2}$} (11.west); 
\path[->,font=\scriptsize,>=angle 90]
(01.south) edge node[right] {$ $} (00.north);
\path[->,font=\scriptsize,>=angle 90]
(11.south) edge node[right] {$ $} (10.north);
\end{tikzpicture}
\end{equation}
is commutative.
\end{lem}

\begin{proof}
This is a straightforward diagram chase, left to the reader.
\end{proof}

\sssec{}

Now consider the case of $E=\O_X$, the trivial line bundle on $X$. Obviously, $\fL_{\O_X} \simeq \pt$. However, on the other hand, the prestack
$$
\fM_{\Ran} 
:=  
\fL_{\O_X,\Ran}
$$
is interesting. Indeed, the previous discussion shows that $\fM_{\Ran}$ is endowed with a structure of commutative monoid, which comes from the system of maps
$$
\O_{X_S}(d D_{\ul x})
\otimes
\O_{X_S}(e D_{\ul y})
\to
\O_{X_S}((d+e) D_{\ul x \cup \ul y}).
$$
Specifically, for $d, e \gg 0$, the above arrow induces a map 
$$
\fM_{X^I}^{\leq d} \times \fM_{X^J}^{\leq e}  \to  \fM_{X^{I \sqcup J}}^{\leq d+e}
$$
of indschemes.
In turn, in the colimit, these maps induce the multiplication $\fM_{\Ran} \times \fM_{\Ran} \to \fM_{\Ran}$. Commutativity is obvious, while associativity is a simple matter.

\begin{rem}
Informally, $\fM_{\Ran}$ classifies pairs $(\ul x, [f])$, with $[f] \in \bbP \bigt{  H^0(X \smallsetminus D_{\ul x}, \O) }$. The monoid structure is given by union on $\Ran$ and by multiplication of functions.
This monoid does not possess a unit element, since $\Ran$ does not either.
\end{rem}

\sssec{}

Similarly, the system of arrows
$$
\O_{X_S}(d D_{\ul x})
\otimes
\restr E{X_S}(e D_{\ul y})
\to
\restr E {X_S} ((d+e) D_{\ul x \cup \ul y})
$$
yields for $d,e \gg 0$ a compatible collection of maps
$$
\fM_{X^I }^{\leq d} \times \fL_{E,X^J }^{\leq e}  \to  \fL_{E,X^{I \sqcup J}}^{\leq d+e}.
$$
In the colimit, these induce an action of $\fM_{\Ran}$ on $\fL_{E,\Ran}$.

\ssec{Taking the quotient}

In the previous sections, we constructed a map $F_{E}: \fL_{E,\Ran} \to \fL_E$ and an action of the monoid $\fM_{\Ran}$ on $\fL_{E, \Ran}$.
In the present section, we combine these two ingredients and show that $\fL_E$ is isomorphic to the fppf sheafification of the quotient $\fL_{E,\Ran}/\fM_{\Ran}$.

\begin{lem}
The map $F_{E}: \fL_{E,\Ran} \to \fL_E$ descends to a map
$$
\fL_{E,\Ran} / \fM_{\Ran} \to \fL_E.
$$
\end{lem}

\begin{proof}
We just need to prove that, for any $d, e \gg 0$, the maps
$$
\fM_{X^I}^{\leq d} \times \fL_{E,X^J}^{\leq e} \xto{\mathit{action}} 
\fL_{E, X^{I \sqcup J}}^{\leq d+e} 
\xto{F_E} 
\fL_E
$$
$$
\fM_{X^I }^{\leq d} \times \fL_{E,X^J}^{\leq e}
\xto{\mathit{projection}}
 \fL_{E,X^J}^{\leq e}
 \xto{F_{E}}
\fL_E
$$
coincide. 
Noting that $\fL_{\O_X} \simeq \pt$, this is an immediate consequence of Lemma \ref{lem:comm-diagram-of-projectivizations}.
\end{proof}

We now claim:

\begin{thm} \label{thm:quotient L/M}
The map 
\begin{equation} \label{eqn:quotient-map}
\phi:
\fL_{E,\Ran} / \fM_{\Ran} \to \fL_E
\end{equation}
is an isomorphism, provided that the quotient is sheafified in the fppf topology.
\end{thm}

The rest of this section (until Section \ref{ssec:Ran-marked}) is devoted to proving the above theorem. The proof is slightly more complicated than expected: this is due to the fact that $\fM_{\Ran}$ is not a group object in $\Ran$-indschemes, but only a non-unital monoid.

\sssec{}

The non-unital monoid $(\Ran, \cup)$ acts on $\fM_{\Ran}$ and $\fL_{E,\Ran}$.
We denote by
$$
\fM_{\indep} := \fM_{\Ran}/\Ran,
\hspace{.4cm}
\fL_{E,\indep} := \fL_{E,\Ran}/\Ran
$$
the quotients (sheafified in the fppf topology).

Observe that $\fM_{\indep}$ is a monoid object in $\PreStk$ and it acts on $\fL_{E,\indep}$. We can thus consider the quotient $\fL_{E,\indep}/\fM_{\indep}$, which by construction comes with a map
$$
\chi:
\fL_{E,\Ran} / \fM_{\Ran} 
\longto
\fL_{E,\indep} / \fM_{\indep}. 
$$
\begin{lem}
The map $\chi$ is an isomorphism.
\end{lem}

\begin{proof}
In view of
$$
\fL_{E,\Ran} / \fM_{\Ran} 
\simeq
\left( 
\fL_{E,\Ran} \times^{\fM_{\Ran}} \fM_{\indep}
\right) 
/ \fM_{\indep},
$$
we just need to prove that the natural arrow
$$
\fL_{E,\Ran} \times^{\fM_{\Ran}} \fM_{\indep}
\longto
\fL_{E,\indep}
$$
is an isomorphism. We have
\begin{eqnarray}
\nonumber
\fL_{E,\Ran} \times^{\fM_{\Ran}} \fM_{\indep}
& \simeq &
\fL_{E,\Ran} \times^{\fM_{\Ran}} 
\left( \fM_{\Ran}/\Ran \right)
\\
\nonumber
& \simeq &
\left( 
\fL_{E,\Ran} \times^{\fM_{\Ran}} 
\fM_{\Ran}  \right) / \Ran
\\
\nonumber
& \simeq &
\fL_{E,\Ran} / \Ran
=: \fL_{E,\indep}
\end{eqnarray}
as desired.
\end{proof}

\sssec{}

It is easy to see that \eqref{eqn:quotient-map} factors as
$$
\fL_{E,\Ran} / \fM_{\Ran} 
\xto{\;\; \chi \;\;}
\fL_{E,\indep} / \fM_{\indep}
\xto{\psi}
\fL_E,
$$
so it remains to prove that $\psi$ is an isomorphism. Since $F_E: \fL_{E,\Ran} \to \fL_E$ is a surjection in the fppf topology, the same is true for the induced map $\fL_{E, \indep} \to \fL_E$.
It follows that $\fL_E$ equals the Zariski sheafification of the colimit of the Cech complex 
\begin{equation}
\label{eqn:Cech-complex}
\cdots 
\fL_{E, \indep} \utimes_{\fL_E} \fL_{E, \indep} \utimes_{\fL_E} \fL_{E, \indep}
\rrr
\fL_{E, \indep} \utimes_{\fL_E} \fL_{E, \indep}
\rr
\fL_{E, \indep}.
\end{equation}

\sssec{}

Our next goal is to describe this Cech complex more explicitly: we wish to see that it is isomorphic to the Cech complex computing the quotient $\fL_{E,\indep}/\fM_{\indep}$ for the $\fM_{\indep}$-action on $\fL_{E,\indep}$.
Concretely, we need to show that 
$$
\fL_{E,\indep} \utimes_{\fL_E} \fL_{E,\indep}
\simeq
\fM_{\indep} \times \fL_{E,\indep}
$$
in such a way that the two projections $\fL_{E,\indep} \utimes_{\fL_E} \fL_{E,\indep} \rr \fL_{E,\indep}$ correspond to projection and action on the right-hand side.
This is the content of the lemmas below.

\begin{lem}\label{lem:qq}
Denote by $r: \fM_{\Ran} \to \Ran$ the structure map. Let $S$ be an affine scheme and $q_1, q_2 \in \fL_{E,\Ran}(S)$ be such that
$
F_{E}(q_1) \simeq F_{E}(q_2).
$
Then there exist an fppf cover $\wt S \tto S$ and $m \in  \fM_{\Ran}(\wt S)$
such that 
$$
r(m) \cdot (\restr {q_1}{\wt S})  = m  \cdot (\restr {q_2}{\wt S})
$$
in $\fL_{E,\Ran}(\wt S)$.
\end{lem}

\begin{proof}

By definition, $q_1 = [\ul x, \eta_1]$ and $q_2 = [\ul y, \eta_2]$ for some $\ul x, \ul y \in \Ran(S)$. Up to acting on both terms by $\ul x \cup \ul y \in \Ran(S)$ and then renaming, we may assume that $\ul x = \ul y$.
Also, up to replacing $S$ with a Zariski cover, we may assume that the line bundles in question are trivial. Thus, for some $d \gg 0$, we have
$$
\eta_1, \eta_2:
\O_S
\hto
(\pi_S)_*( \restr E {X_S}(d D_{\ul x})).
$$
By construction, for $i =1,2$, the element $F_{E}(q_i)$ is represented by a pair $[U_i, \eta'_i]$, where $U_i \Subset U_{\ul x} \Subset X_S$ and
$$
\eta'_i: \O_{U_i} \hto \restr E {U_i}
$$
is the line sub-bundle obtained from $\eta_i$ as in Section \ref{ssec:FE-construction}. Now, the assumption $F_\infty(\ell_1) = F_{\infty}(\ell_2)$ implies that $ \restr {\eta_1} U \simeq  \restr {\eta_2} U $ on some subdomain $U \Subset U_1 \cap U_2$.
Thus, there exists an invertible function $h \in H^0(U,\O_U)$ such that
\begin{equation} \label{eqn:eta = eta times h}
\eta'_1|_U =  \eta'_2|_U \circ (h \cdot -)
:  \O_U \to \restr E U.
\end{equation}
We now argue as in Proposition \ref{prop:FE-epi}: by further shrinking $U$ and by pulling back along an appropriate fppf cover $\wt S \to S$ and  we may assume that $U = U_{\ul z}$ for some $\ul z \in \Ran(\wt S)$ containing $\ul x$.
Thus, the pullback of $h$ to $\wt S$ is represented by a section $\wt h \in H^0(X_{\wt S}, \O_{X_{\wt S}}(e D_{\ul z}))$ for some $e \gg 0$. 
By adjunction, $\wt h$ determines a map
$$
\wt h':
\O_{\wt S} \to (\pi_{\wt S})_*
\bigt{ \O_{X_{\wt S}}(e D_{\ul z})
},
$$
which is an injection of coherent sheaves with locally free cokernel.
We obtain an element $m: = [\ul z,\wt  h'] \in \fM_{\Ran}(\wt S)$ and claim that
$$
\ul z \cdot \restr {q_1} {\wt S} = m \cdot \restr {q_2}{\wt S}.
$$
This is obvious: by \eqref{eqn:eta = eta times h} the equality is valid on $U_{\ul z}$, hence it is valid on $X_{\wt S}$ by (schematic) density. 
\end{proof}

\begin{lem}
Suppose that $m_1 \cdot q = m_2 \cdot q$ for some
$m_1, m_2 \in \fM_{\Ran}(S)$ and $q \in \fL_{E, \Ran}(S)$. Then there are $r_1, r_2 \in \Ran(S)$ such that 
$
r_1 \cdot m_1 = r_2 \cdot m_2.
$
\end{lem}

\begin{proof}
Keeping our usual notation, we set
$$
m_1 = (\ul x, \zeta_1), 
\hspace{.4cm}
m_2 = (\ul y, \zeta_2),
\hspace{.4cm}
q = (\ul z, \eta).
$$
Up to acting on all terms by $\ul x \cup \ul y \cup \ul z \in \Ran(S)$ and then renaming, we may assume that $\ul x = \ul y =\ul z$. With this modification, we claim that $m_1 = m_2$.
There are line bundles $L_1, L_2, L_3$ on $S$ and  $d \gg 0$ large enough such that $m_i$ (for $i=1,2$) corresponds to
$$
\zeta_i \in 
\Hom
_{\Coh(X_S)^{\heart}}
\Bigt{
\restr{L_i}{X_S},\O_{X_S} (d D_{\ul x}) 
},
$$
while $q$ corresponds to
$$
\eta \in \Hom_{\Coh(X_S)^{\heart}}
\Bigt{
\restr{L_3}{X_S} ,\restr E {X_S} (d D_{\ul x}) 
}.
$$
In turn, the elements $m_i \cdot q$ (for $i=1,2$) correspond to
$$
\zeta_i \cdot \eta \in \Hom_{\Coh(X_S)^{\heart}}
\Bigt{
\restr{(L_i \otimes L_3)}{X_S},\restr E {X_S} (2d D_{\ul x}) 
}.
$$
The condition that $m_1 \cdot q = m_2 \cdot q$ implies the existence of an isomorphism $\theta: L_1 \xto\simeq L_2$ such that 
$$
(\zeta_1 \circ \widehat \theta) \cdot \eta  = \zeta_2 \cdot \eta,
$$ 
where $\widehat \theta: \restr{(L_1 \otimes L_3)}{X_S} \to \restr{(L_2 \otimes L_3)}{X_S} $ is the isomorphism induced by $\theta$.
Now, by assumption, $\eta$ is injective as a map of coherent sheaves: thisimplies that $\zeta_1 \circ \widehat \theta = \zeta_2$ on $X_S \smallsetminus D_{\ul x}$, and consequently $\zeta_1 \circ \widehat \theta = \zeta_2$ on the entire $X_S$. By definition, this means that $m_1 = m_2$, as required.
\end{proof}

\ssec{Ran space with marked points} \label{ssec:Ran-marked}

This subsection plays a technical role and can be skipped by the reader until needed. 

\sssec{}

Let $D_0 = \{z_1, \ldots, z_m\}$ be a finite set of $k$-points of $X$, fixed once and for all.
Following \cite[Section 3.6.6]{contract}, we will define a space $\mRan$ that parametrizes finite subsets of $X$ that contain $D_0$. 
Namely, for $I \in \fSet$, we denote by $X^{I, \ast}$ the closed subscheme of $X^I$ consisting of those $I$-tuples that contain $D_0$ as a subset. We then set
$$
\mRan := \colim_{I \in \fSet^\op} X^{I,\ast}.
$$
The usual proof shows that $\mRan$ is both $\O$-contractible and homologically contractible.

\sssec{}

We leave it to the meticulous reader to check that the constructions of the previous part of Section \ref{sec:contract-LE} can be rerun with $\mRan$ in the place of $\Ran$.
In particular, we have $\mRan$-indschemes $\fL_{E, \mRan}$ and $\fM_{\mRan}$, an action of the latter on the former, and an isomorphism
$$
\fL_{E, \mRan}/ \fM_{\mRan} \xto{\;\; \simeq \;\;} \fL_E,
$$
with the quotient sheafified in the fppf topology.
Furthermore, $\fL_{E, \mRan}$ and $\fM_{\mRan}$ are both $\O$-contractible.

\ssec{The contractibility of $\fL$}

Here we provide a new proof of the homological contractibility of $\fL_E = \GSect(X,\bbP E)$. We do this by combining two other results that are of interest in their own right.  The first result is the $\O$-contractibility of $\fL_E$.
The second one is the fact that the forgetful functor
$$
\oblv: 
\Dmod(\fL_E)
\longto
\QCoh(\fL_E)
$$
is fully faithful. 

\sssec{}

Recall that, in the case of pseudo-proper pseudo-schemes like $\fL_E$, the $\O$-contractibility and the $\ICoh$-contractibility are equivalent.

\begin{thm} \label{thm:O-contract}
For any vector bundle $E \to X$, the pseudo-proper pseudo-scheme $\fL_E$ is $\ICoh$-contractible (and thus $\O$-contractible).
\end{thm}

\begin{proof}
Since $\ICoh$ satisfies fppf descent, the quotient presentation of $\fL_E \simeq \fL_{E, \Ran}/\fM_{\Ran}$ proven in Theorem \ref{thm:quotient L/M} yields
$$
\ICoh(\fL_E)
\simeq
\colim_{[n] \in \Delta^\op}
\Bigt{
\ICoh(\fM_{\Ran})^{\otimes n} \otimes
\ICoh(\fL_{E,\Ran})
},
$$
where the colimit is taken with respect to the pushforward functors. This formally implies that
$$
H_*(\fL_E,\O)
\simeq
\colim_{[n] \in {\Delta}^\op}
\Bigt{
H_*(\fM_{\Ran}, \O)^{\otimes n} \otimes
H_*(\fL_{E, \Ran}, \O)
}.
$$
By Proposition \ref{prop:fL-E-Ran contractible}, the spaces $\fL_{E,\Ran}$ and $\fM_{\Ran}$ are both $\ICoh$-contractible. If follows that 
$$
H_*(\fL_E,\O)
\simeq
\colim_{[n] \in \Delta^{\op}} k
\simeq
k,
$$
the latter isomorphism following from the fact that $\Delta^{\op}$ is a contractible category.
\end{proof}

\begin{lem}
The standard comparison functor $\Upsilon_{\fL_E}: \QCoh(\fL_E) \to \ICoh(\fL_E)$ is an equivalence.
 \end{lem}
 
\begin{proof}
Since $\QCoh$ and $\ICoh$ satisfy fppf descent, we just need to show that $\Upsilon_{\fL_{E,\Ran}}$ and $\Upsilon_{\fM_{\Ran}}$ are equivalences. Since $\fM_{\Ran}$ is a particular instance of ${\fL_{E,\Ran}}$, we just focus on the latter. In this case, the assertion follows from the fact that each $\fL_{E, X^I}$ is formally smooth.
\end{proof}

\sssec{}

Next, we prove the following.

\begin{thm} \label{thm:fully-faith}
For any vector bundle $E \to X$, the forgetful functor
$$
\oblv: 
\Dmod(\fL_E)
\longto
\QCoh(\fL_E)
$$
is fully faithful.
\end{thm}

\begin{proof}

We proceed in steps.

\sssec*{Step 0}
Let us begin with some preparation.
In view of the lemma above, we can post-compose with $\Upsilon_{\fL_E}$ and prove that the forgetful functor
$$
\oblv: 
\Dmod(\fL_E)
\longto
\ICoh(\fL_E)
$$
is fully faithful. This passage to $\ICoh$ is convenient for a technical reason: we will need to work with a formal completion on which we prefer to consider ind-coherent sheaves.

Furthermore, by Remark \ref{rem:can-trivialize-E}, we can assume that $E$ is trivial, so that $\fL_E \simeq \GMaps(X, \bbP^n)$ for $n = \mathrm{rank}(E)-1$. To reduce clutter, we set $\fL := \GMaps(X, \bbP^n)$.

\sssec*{Step 1}

The functor $\oblv$ is, by definition, given by the ind-coherent pullback along the tautological map $q: \fL \to \fL_\dR$. Since this map is a nil-isomorphism, it suffices to prove that the canonical arrow
$$
q_*^{\ICoh}(\omega_{\fL}^{\ICoh})
\longto
\omega_{\fL}
$$
is an isomorphism in $\Dmod(\fL)$. We can check this on $k'$-fibers for extensions $k' \supseteq k$.
Applying base-change, it remains to prove that, for any $\ell \in \fL(k')$, the map 
$$
\Gamma^{\ICoh}_{k'}
\left(
(\fL \otimes_k k')^\wedge_\ell, \omega
\right ) 
 \to k'
$$ 
is an isomorphism. 
Now, by renaming $\fL \otimes_k k'$ by $\fL$ and $k'$ by $k$, we may assume that $\ell \in \fL(k)$ and thus we just need to prove that
$$
\Gamma^{\ICoh}\left(  \fL^\wedge_\ell, \omega \right) 
\simeq 
k.
$$

\sssec*{Step 2}

Now consider the action of $\GMaps(X, GL_{n+1})$ on $\fL = \GMaps(X, \bbP^n)$. This action is evidently transitive at the level of $k$-points, so we can assume that $\ell$ is the $k$-point of $\fL$ corresponding to the first coordinate.

Now, we use the presentation $\fL = \fL_{\Ran}/\fM_{\Ran}$ of Theorem \ref{thm:quotient L/M}.
Since $E = \O_X^{\oplus n+1}$, the canonical injection $\O_X \hto \O_{X}^{\oplus n+ 1}$ in the first coordinate induces a closed embedding
$$
e_1: \fM_{\Ran} \hto \fL_{\Ran}.
$$ 
By construction, $e_1$ yields the $k$-point $\ell$ at the level of $\fM_{\Ran}$-quotients.
It follows that
$$
\fL^\wedge_\ell
\simeq
(\fL_{\Ran})^\wedge_{\fM_{\Ran}}/\fM_{\Ran},
$$
with the formal completion of course taken with respect to $e_1$.
Since $\fM_\infty$ is $\ICoh$-contractible, it suffices to prove that 
$$
\Gamma^{\ICoh}
\bigt{
(\fL_{\Ran})^\wedge_{\fM_{\Ran}}, \omega
} 
\simeq k.
$$

\sssec*{Step 3}

For $I \in \fSet$, consider the structure map $p_I: (\fL_{X^I})^\wedge_{\fM_{X^I}} \to X^I$. It suffices to prove that the natural arrow
$$
(p_I)_*^{\ICoh}
\Bigt{
\omega_{(\fL_{X^I})^\wedge_{\fM_{X^I}}}
}
\to
\omega_{X^I}
$$
is an isomorphism. Indeed, this way, we would have
$$
\Gamma^{\ICoh}
\bigt{
(\fL_{\Ran})^\wedge_{\fM_{\Ran}}, \omega
} 
\xto{\simeq}
\Gamma^{\ICoh}(\Ran,\omega)
$$
and we could conclude thanks to the $\ICoh$-contractibility of $\Ran$.

\sssec*{Step 4}

Now recall that $\fL_{X^I}$ and $\fM_{X^I}$ are both indschemes in a natural way: 
$$
\fL_{X^I} \simeq \colim_{d \gg 0} \fL_{X^I}^{\leq d}
$$
and similarly for $\fM_{X^I}$. Hence, it suffices to find a divergent sequence $\{N_d\}_{d }$ of natural numbers such that 
\begin{equation} \label{eqn:coker-fin-dim}
\cone
\left(
(p_I^{\leq d})_*^{\ICoh}
\Bigt{
\omega_{(\fL_{X^I}^{\leq d})^\wedge_{\fM_{X^I}^{\leq d}}}
}
\to
\omega_{X^I}
\right)
\in \ICoh(X^I)^{\leq - N_d}
\end{equation}
for all $d \gg 0$.

\sssec*{Step 5}

Let choose $d_0$ in such a way that, for all $d \geq d_0$, the schemes $\fL_{X^I}^{\leq d}$ and $\fM_{X^I}^{\leq d}$ are both projective bundles over $X^I$, see Lemma \ref{lem:locally free}. Since $X^I$ is quasi-compact, it suffices to find a suitable $N_d$ separately on each member of a Zariski cover of $X^I$. Choose an open subset $U \subset X^I$ on which projective bundles are trivial. Computing the  dimension of the fibers by Riemann-Roch, we obtain that
$$
\restr
{
\fL_{X^I}^{\leq d}}
{U}
\simeq
U \times \bbP^{L_d}
$$
$$
\restr
{
\fM_{X^I}^{\leq d}}
{U}
\simeq
U \times \bbP^{M_d }
$$
with
$$
L_d: = (|I| d +1 -g) (n+1) - 1
= ((n+1) |I|) \cdot d + O(1)
$$
$$
M_d := ( |I| d +1 -g) - 1 
=
|I| \cdot d + O(1).
$$
These isomorphisms obviously can be arranged so that the closed embedding in question is the standard linear embedding $\bbP^{M_d} \hto \bbP^{L_d}$.

\sssec*{Step 6}

Then the object appearing in \eqref{eqn:coker-fin-dim}, restricted to $U$, is isomorphic to
$$
\cone
\left(
\Gamma^{\ICoh}
\Bigt{
(\bbP^{L_d})^\wedge_{M_d}, \omega
} 
\to k
\right)
\otimes \omega_U.
$$
It suffices to show that the cone in question is supported cohomologically in degrees $\leq - M_d$. This is a computation of general interest, which is recorded in the lemma below.
\end{proof}

\begin{lem}
For any pair of natural numbers $a \geq b$, consider the standard closed embedding $\bbP^b \hto \bbP^a$. The cone of the map
$$
\Gamma^{\ICoh}
\Bigt{
(\bbP^a)^\wedge_{\bbP^b}, \omega
} 
\to k
$$
belongs to $\Vect^{\leq -b}$.
\end{lem}

\begin{proof}
The infinitesimal neighbourhoods of the closed subscheme $\bbP^b \subseteq \bbP^a$ yield a non-negative filtration on $\Gamma^{\ICoh}
\Bigt{
(\bbP^a)^\wedge_{\bbP^b}, \omega
} $ with $n$-associated graded given by
$$
\Gamma^{\ICoh}
(
\bbP^b, \Sym_{\bbP^b}^n( \bbT_{rel}[1])
).
$$
By Serre duality, the latter is equal to the dual of 
$$
H^*(\bbP^b, \Sym_{\bbP^b}^n( \CN^\vee)),
$$
where $\CN^\vee$ is the conormal bundle. 
It suffices to prove that, for all $n \geq 1$, this chain complex has cohomologies only in degrees $\geq b$. This is clear, since $\CN^\vee$ is a direct sum of copies of $\O_{\bbP^b}(-1)$.
\end{proof}

\sssec{}

The following corollary is well-known, see \cite{Barlev}, but our proof is new.

\begin{cor}
For any vector bundle $E \to X$, the pseudo-scheme $\fL_E$ is homologically contractible.
\end{cor}

\begin{proof}
Theorem \ref{thm:fully-faith} shows that 
$$
H^*(\fL_E)
:=
\CHom_{\Dmod(\fL)}(\omega_{\fL_E}, \omega_{\fL_E})
\xto{\simeq}
\CHom_{\QCoh(\fL)}(\O_{\fL_E}, \O_{\fL_E})
\simeq
H^*(\fL_E, \O).
$$
The latter is isomorphic to $k$ by the $\O$-contractibility of $\fL_E$ of Theorem \ref{thm:O-contract}.
\end{proof}

\sec{Proof of the main theorem} \label{sec:Oka}

In the previous section, we established the homological contractibility of the space $\Opbargen_{G,\sigma}$ for general linear groups. It remains to complete the proof of Theorem \ref{mainthm:contract-fiber} for orthogonal and symplectic groups. 
To do this, we use a general method, described below, that lets us reduce the question to analyzing an explicit closed subspace of $\fL_{E, \Ran}$. We then conclude the proof by performing a t-structure estimate of weak-Lefschetz type.

\ssec{Contractibility of subspaces of $\fL_E$}  \label{ssec:subspaces}

Let us provide a general method to test the contractibility of subspaces of $\fL_E$. 

\sssec{}

Let $j: \fB \hto \fL_E$ be a locally closed embedding. We set
$$
\fB_{\Ran} 
:=
\fB \times_{\fL_E} \fL_{E, \Ran}.
$$
The action of the monoid $\fM_{\Ran}$ on $\fL_{E, \Ran}$ obviously preserves $\fB_{\Ran}$, hence we obtain 
$$
\fB 
\simeq
\fB_{\Ran}  / \fM_{\Ran} .
$$

\begin{prop}\label{prop:contractibility-at-Ran-level}
If the map $j_{\Ran}: \fB_{\Ran}  \to \fL_{E,\Ran} $ induces an isomorphism
$$
H_*(\fB_{\Ran}) \xto{\simeq} H_*(\fL_{E,\Ran} ),
$$
then $\fB$ is homologically contractible.
\end{prop}

\begin{proof}
The isomorphism
$\fB 
\simeq
\fB_{\Ran}  / \fM_{\Ran} 
$
yields an equivalence
$$
\Dmod(\fB)
\simeq
\colim_{[n] \in \Delta^{\op}}
\bigt{
\Dmod(\fM_{\Ran} )^{\otimes n} 
\otimes 
\Dmod(\fB_{\Ran} )
},
$$
and consequently an isomorphism
$$
H_*(\fB)
\simeq
\colim_{[n] \in \Delta^{\op}}
\bigt{
H_*(\fM_{\Ran} )^{\otimes n} \otimes H_*(\fB_{\Ran} )
}.
$$
Then the isomorphism $H_*(\fB_{\Ran}) \simeq H_*(\fL_{E,\Ran} )$, which is necessary compatible with the $H_*(\fM_{\Ran})$-action, yields an isomorphism $H_*(\fB) \simeq H_*(\fL_E)$. We conclude by the homological contractibility of $H_*(\fL_E)$.
\end{proof}

\begin{rem}[This can be skipped by the reader]

The same idea yields the following result as well.
Let $i: \fB \hto \fL_E$ be a closed embedding, and suppose that the induced map $i_{\Ran}: \fB_{\Ran}  \hto \fL_{E,\Ran} $ yields an isomorphism
$$
H_*(\fB_{\Ran}, \O) \simeq H_*(\fL_{E,\Ran},\O).
$$
Recall that these vector spaces are pre-dual to the more familiar coherent cohomologies. Then $\fB$ is $\O$-contractible.
Indeed, the isomorphism
$\fB 
\simeq
\fB_{\Ran}  / \fM_{\Ran},
$
together with the $\ICoh$-contractibility of $\fM_{\Ran}$, yields an isomorphism $H_*(\fB_{\Ran}, \O) \simeq H_*(\fB, \O)$. On the other hand, we already know that $\fL_{E, \Ran}$ is $\ICoh$-contractible.

\end{rem}

\sssec{} \label{sssec:favourable situations}

Let us return to disucussing the homological contractibility of $\fB$. In favourable situations, we can hope to prove a statement stronger than the above proposition, as follows. Let $\pi_{\Ran}: \fL_{E,\Ran} \to \Ran$ denote the structure map. Retaining the above notation, suppose that the object
\begin{equation} \label{eqn:vanishing on Ran}
(\pi_{\Ran})_!
\Bigt{
\cone
\bigt{
(j_{\Ran})_!(\omega_{\fB_{\Ran}})
\to
\omega_{\fL_{E,\Ran}}
}
}
\in \Dmod(\Ran)
\end{equation}
vanishes. Then the assumption of Proposition \ref{prop:contractibility-at-Ran-level} is satisfied and thus $\fB$ is homologically contractible.

\sssec{}

Now, the vanishing of \eqref{eqn:vanishing on Ran} can be checked on field-valued $!$-fibers. If $j_{\Ran}$ is a closed embedding, we deduce the following result by proper base-change.

\begin{cor}
Let $i: \fB \hto \fL_E$ be a closed embedding. Suppose that, for any $k' \supseteq k$ and any $\ul x: \Spec(k') \to \Ran$, the induced arrow
$$
H_*(\fB_{\ul x}) \to H_*(\fL_{E, \ul x})
$$
of $k'$-vector spaces is an isomorphism. Then $\fB$ is homologically contractible.
(Here $\fB_{\ul x}$ is the pseudo-scheme over $k'$ defined as $\fB_{\Ran} \times_{\Ran} \Spec(k')$, and similarly for $\fL_{E,\ul x}$.)
\end{cor}

\sssec{}

The statement of the above corollary can be approached by finite dimensional geometry.
Namely, the map $i_{\ul x}: \fB_{\ul x}  \hto \fL_{E,\ul x} $ arises as a colimit of closed embeddings
$$
i^{\leq d}_{\ul x}: 
 \fB_{\ul x}^{\leq d}  \hto \fL_{E,\ul x}^{\leq d} 
$$
where
$$
\fB_{\ul x}^{\leq d} 
:=
\fB_{\ul x}
\times_{\fL_{E, \ul x}}
\fL_{E,\ul x}^{\leq d}.
$$
Note that the $k'$-scheme $\fL_{E,\ul x}^{\leq d}$ is a finite dimensional projective space, the projectivization of the (finite dimensional) $k'$-vector space 
$H^0 (X_{k'}, \restr E{X_{k'}} (d D_{ \ul x}))$.

\begin{cor} \label{cor:contractib-x-d}
Let $i: \fB \hto \fL_E$ be a closed embedding and $\ul x$ a $k'$-point of $\Ran$.
Suppose there exists a divergent sequence $\{N_d\}_{d \geq 0}$ of natural numbers such that, for any $d \geq 0$, the cone of the map
$$
H_*(\fB_{\ul x}^{\leq d})
\to 
H_*(\fL_{E,\ul x}^{\leq d})
$$
belongs to $\Vect_{k'}^{\leq - N_d}$.
Then $\fB$ is homologically contractible.
\end{cor}

\begin{proof}
We will deduce this claim from the previous corollary. By switching colimits, we see that 
$$
\cone
\bigt{
H_*(\fB_{\ul x}) \to H_*(\fL_{E,\ul x} )
}
\simeq
\colim_{d \gg 0}
\Bigt{
\cone
\bigt{
H_*(\fB^{\leq d}_{\ul x}) \to H_*(\fL_{E,\ul x}^{\leq d})
}
}.
$$
Then the assumption of our corollary implies that this colimit belongs to $\Vect_{k'}^{\leq - N_d}$ for any $d \gg 0$. Since the t-structure on the DG category $\Vect_{k'}$ is left-complete and $\lim_{d \to \infty} N_d = \infty$, the colimit must vanish.
\end{proof}

\begin{rem} \label{rem:marked-contract}
The hypotheses of the above corollary can be slightly relaxed. Namely, take any finite subset $ D_0 \subset X(k)$ and consider $\mRan$, the marked Ran space of $(X,D_0)$. Then, in view of Section \ref{ssec:Ran-marked}, the assertion of the corollary holds with $\Ran$ replaced by $\mRan$.
\end{rem}

\ssec{Oka's theorem (or t-structure estimates)}

We now apply the above theory to $\Opbargen_{G,\sigma}$. For definiteness, we focus only on the case $G=Sp_{2n}$, the situation for the orthogonal groups is similar and left to the reader.

\sssec{} \label{sssec:our-marking}

Thus, let $G = Sp_{2n}$ and $\sigma = (E,\nabla, \omega)$ a $G$-local system on $X$. 
As in Section \ref{ssec:type-Cn}, we can puncture $X$ by removing a finite subset $D_0 = \{z_1, \ldots, z_m\}  \subset X(k)$ to obtain an affine open subvariety $X' \subseteq X$ such that $\restr E {X'}$ and $T_{X'}$ are both trivialized.

The set $D_0$ is the marking we use when involving the marked Ran space $\mRan$ in the proof of Proposition \ref{prop:estimate-opers} below.

\sssec{}

Before getting to the proof, we need to fix some more notation. Given an integer $d\geq 0$, we denote
$$
k^{d}[X'] 
:=
H^0( X ,\O_{ X}(d D_0)).
$$
Note that $k^d[X']$ is a finite-dimensional $k$-vector space: by Riemann-Roch, for $d \gg 0$, we have 
$
\dim_k k^d[X] =   d|D| - g +1. 
$
We also have
$$
k[X'] := H^0(X', \O_{X'}) 
\simeq 
\colim_{d \gg 0} k^d[X']. 
$$

\sssec{} \label{sssec:vector field}

Let us also fix a nowhere vanishing vector field $\nu$ on $X'$ that trivializes the tangent bundle $T_{X'}$.
We view $\nu$ as a derivation $\nu: k[X'] \to k[X']$. Now fix an open subscheme $U \subseteq X'$ and consider $\restr \nu U$. It is easy to see that there exists a constant $d_{\nu} \geq 0$ (depending only on $U$) such that, for every $d\geq 0$, we have
$$
\nu(k^d[U]) \subset k^{d+d_{\nu}}[U].
$$

\sssec{}

In Section \ref{ssec:type-Cn}, we constructed a pseudo-scheme $\Y_{G, \sigma}$ and a $\xi: \Opbargen_{G,\sigma} \to \Y_{G, \sigma}$ inducing an isomorphism in homology. The space $\Y_{G,\sigma}$ comes with a closed embedding 
$$
i:
\Y_{G,\sigma}
\hto
\GMaps(X,\bbP V),
$$
where $V$ is the standard $G$-representation.
Set $\fL := \GMaps(X,\bbP V)$.
Our goal is to show that the induced map $H_*(\Y_{G,\sigma}) \to
H_*(\GMaps(X,\bbP V))$ is an isomorphism.
We will prove this by employing the general method of Section \ref{ssec:subspaces}.

\sssec{}

Namely, up to fppf sheafification, we have:
$$
\fL
\simeq
\fL_{\Ran} /\fM_{\Ran},
$$
$$
\Y_{G,\sigma}
\simeq
\fY_{\Ran} /\fM_{\Ran},
$$
where
$$
\fY_{\Ran}
:=
\fL_{\Ran} \utimes_{\fL} \Opbargen_{G,\sigma}.
$$
In fact, we will need to consider the above formulas for the marked Ran space relative to our fixed marking $D_0 \subset X(k)$ of Section \ref{sssec:our-marking}.

Then, in view of Corollary \ref{cor:contractib-x-d} and Remark \ref{rem:marked-contract}, the proof of Theorem \ref{mainthm:contract-fiber} boils down to the cohomological estimate below.

\begin{prop} \label{prop:estimate-opers}
Let $k' \supseteq k$ a field extension and $k'$-point $\ul x$ of $\mRan$. There exists a divergent sequence $\{N_d\}_{d \geq 0}$ of natural numbers such that, for any $d \gg 0$, the cone of the map
$$
H_*(\fY_{\ul x}^{\leq d})
\to 
H_*(\fL_{\ul x}^{\leq d})
$$
belongs to $\Vect_{k'}^{\leq - N_d}$.
\end{prop}

\begin{proof}
First of all, by renaming $k$ with $k'$, we may assume that $\ul x \in \mRan(k)$. By construction, the support of $D:= D_{\ul x}$ contains $D_0$.
Letting $U := X \smallsetminus D$, we have an inclusion $U \subseteq X' := X \smallsetminus D_0$ in the opposite direction.

\sssec*{Step 1}

By definition,  $\fL_{\ul x}^{\leq d}
\simeq \bbP(H^0(X, \O_X^{\oplus 2n} (dD))$.
Thus, for $d \gg 0$, Riemann-Roch yields
$$
\fL_{\ul x}^{\leq d}
\simeq 
\bbP^{( 2n|D|) d + c_0 },
$$
where $c_0 = 2n(1-g)-1$ is a constant (obviously) independent of $d$.

\sssec*{Step 2}

Since $U \subseteq X'$ and $\restr E {X'}$ has been trivialized, we can write the alternating form
$$
\restr \omega U :  \Lambda^2 ( \O_U^{\oplus 2n} ) \to \O_U
$$
as a matrix
$
M_\omega \in M_{2n}(k[U])
$.
Similarly, the differential operator
$$
\restr{\nabla_\nu} U : \O_U^{\oplus 2n}
 \to 
 \O_U^{\oplus 2n}
$$
can be written as $\nu + A$ for some matrix $A \in M_{2n}(k[U])$.
In particular, there exists a constant $c_1 \in \bbN$ be such that all the entries of $M_\omega$ and of $A$ lie in $k[U]^{c_1} := H^0(X, \O(c_1  D))$.

\sssec*{Step 3}

Next, by Section \ref{sssec:vector field}, we see that there is some constant $d_\nu \in \bbN$ (independent of $d$) such that the derivation
$
\nu : k[U] \to k[U]
$
restricts, for every $e \geq 0$, to a map
$$
\nu:  k[U]^{e} \to k[U]^{e + d_\nu}.
$$
Combining these estimates, we can find a constant $c_2 \in \bbN$ such that, for any $d \geq 0$ and any section $v\in H^0(X,\O_X^{\oplus 2 n}(dD))$, the $(n-1)$ meromorphic functions
\begin{eqnarray}
\nonumber
& &
\omega ( v,\nabla_\nu(v)) \\
\nonumber
& &
\omega ( \nabla_\nu (v),\nabla_\nu^{(2)}(v)) \\
\nonumber
& &
\cdots \\
\nonumber
& &
\omega(\nabla_\nu^{(n-2)} v,\nabla_\nu^{(n-1)} v)
\end{eqnarray}
all lie in $k[U]^{2d+c_2} = H^0(X, \O((2d+c_2)  D))$.

\sssec*{Step 4}

By definition, $\fY_{\ul x}^{\leq d}$ is the closed sub-scheme of $\fL_{\ul x}^{\leq d}$ cut out by the vanishing of the above displayed functions. By Riemann-Roch again, these conditions amount to $2(n-1)|D|d + c$ equations, for some $c$ independent of $d$.

Denote by 
$$
N_d := \dim(\fL_{\ul x}^{\leq d}) - (2(n-1)|D|d + c)
=
2|D|d + 
O(1)
$$
the difference between the dimension of the ambient projective space and the number of equations. 
On the one hand, $N_d$ tends to $\infty$ with $d$. On the other hand, the following result of \virg{weak Lefschetz} type shows that
$$
\cone
\Bigt{
H_*(\fY_{\ul x}^{\leq d}) 
\to H_*(\fL_{\ul x}^{\leq d})
}
\in
\Vect^{\leq - N_d}
$$
as desired.
\end{proof}

\begin{lem}[Oka]
Let $V$ be an algebraic set in $\bbP^N$ defined by $r$ homogeneous polynomials. Then the cone of the map 
$$ 
H_*(V) \to H_*(\mathbb{P}^N)
$$
belongs to $\Vect^{\leq - (N-r)}$.
\end{lem}

\begin{proof}
We will give two proofs, the first one is topological, the second one more algebraic.

\sssec*{The first proof}

Since the definition of $V$ involves only finitely many coefficients, we may assume that $k$ is of finite transcendence degree and thus that it admits an inclusion $k \hto \bbC$. We can then apply the Riemann-Hilbert correspondence and work in the classical analytic topology. In particular, we can now invoke Oka's result, see \cite{Oka} and \cite[Lemma 6.1]{Kato}.

\nc{\wh}{\widehat}

\sssec*{The second proof}

Consider the affine cone $\wh{V} \subseteq \bbA^{N+1}$. Let us work with the \virg{perverse} t-structure on the DG category of D-modules. By \cite[Lemma 6.2.2]{omega}, the object $\omega_{\wh{V}}$ belongs to $\Dmod(\wh V)^{\leq -(N-r)}$ and, by affineness, $\Hren(\wh V) \in \Vect^{-(N -r)}$.

Now, the blow-up correspondence between $\wh V$ and $V$, together with the projective bundle formula, yields the isomorphism
$$
\Hren(\wh V) \oplus \ol H_*(V)
\simeq
H_*(V)[2],
$$
where $\ol H_*$ denotes the reduced homology.
In particular, $\ol H_{-i}(V) \simeq H_{-i-2}(V)$ for all integers $i < N-r$ and the assertion follows.
\end{proof}

\sec{A higher Tsen theorem} \label{sec:Tsen}

As mentioned, our methods in the proof of Theorem \ref{mainthm:contract-fiber} are quite general. As an illustration of this, we apply them to prove a homotopical strengthening of Tsen's theorem about rational points over function fields of curves.

\setcounter{subsection}{1}

\sssec{}
 
Let us first recall the statement of Tsen's theorem.
\begin{thm}[Tsen]
Let $k$ be an algebraically closed field and let $X$ be a smooth curve defined over $k$.
Let $Z \subset \mathbb{P}^n_{k(X)}$ be the zero locus of $r$ homogeneous polynomials $f_i \in k(X)[t_0, \ldots, t_n]$ of degrees $d_1,\dots, d_r$. If
$\sum_{i=1}^r d_i \leq n$, then $Z(k(X)) \neq \emptyset$.
\end{thm}

\sssec{}

In our strengthening of Tsen's theorem, we wish to show that $Z(k(X))$ is not just nonempty but also homologically contractible.
To this end, we will construct a pseudo-scheme $\fZ$ over $k$ (dependent only on $Z$) with $\fZ(k) \simeq Z(k(X))$ and show that $\fZ$ is homologically contractible.

\sssec{}

To define $\fZ$, consider the polynomials $f_1, \ldots, f_r$. For an appropriate large effective divisor $D_0$ on $X$ we can view each $f_i$ as a section 
$$
F_i \in H^0 \bigt{ X \times \bbP^n, \O_X(D_0) \boxtimes \O_{\bbP^n}(d_i)}.
$$
We then define the closed subscheme $\wt Z \subseteq \bbP^n_X$ to be the zero locus of the $F_i$'s. Next, we set
$$
\fZ := \GSect(X, \wt Z).
$$
Note that $\fZ$ depends only on $Z$, and not on $D_0$. Reasoning along the same lines, it is clear that $\fZ(k) = Z(k(X))$.

\begin{thm}
If $\sum d_i \leq n$, then the pseudo-scheme $\fZ:= \GSect(X,\wt Z)$ is homologically contractible.
\end{thm}

\begin{proof}
We can present 
$$
\GSect(X,\CZ) \subseteq \GMaps(X, \bbP^n) =: \fL
$$ 
as a quotient of a $\Ran$-indscheme $\fZ_{\Ran}$ by the monoid $\fM_{\Ran}$.
As in the proof of Proposition \ref{prop:estimate-opers}, we just need to prove the following claim: for $\ul x \in \Ran(k)$, the closed subscheme $\fZ_{\ul x}^{\leq d} \subseteq \fL_{\ul x}^{\leq d}$ is defined as the zero locus of $n_d$ equations with 
$$
N_d := \dim (\fL_{\ul x}^{\leq d}) - n_d
$$
going to $\infty$ with $d$. Since
$$
\dim (\fL_{\ul x}^{\leq d}) = (n+1)|D| \cdot d + O(1)
$$
$$
n_d = (\sum d_i) |D| \cdot d + O(1),
$$
our assertion follows.
\end{proof}

\end{document}